\PassOptionsToPackage{table}{xcolor}
\documentclass[onefignum,onetabnum,a4paper]{siamart190516}
\pdfoutput=1
\usepackage[utf8]{inputenc}
\usepackage[english]{babel}
\usepackage{amsmath, amssymb, gensymb, physics, upgreek, siunitx, bbm}
\usepackage{tikz, pgfplots, graphicx, epstopdf, float, hyperref, xspace}
\usepackage[caption=false]{subfig}
\usepackage{enumitem, verbatim, listings}
\usepackage{mathtools, stmaryrd}
\usepackage{booktabs, array, ragged2e}
\usepackage[export]{adjustbox}
\usepackage{cleveref}

\usetikzlibrary{trees, cd, babel}
\graphicspath{{./figures}}
\definecolor{paired1}{HTML}{a6cee3}
\definecolor{paired2}{HTML}{1f78b4}
\definecolor{paired3}{HTML}{b2df8a}
\definecolor{paired4}{HTML}{33a02c}
\definecolor{paired5}{HTML}{fb9a99}
\definecolor{paired6}{HTML}{e31a1c}
\definecolor{paired7}{HTML}{fdbf6f}
\definecolor{paired8}{HTML}{ff7f00}
\definecolor{paired9}{HTML}{cab2d6}
\definecolor{paired10}{HTML}{6a3d9a}

\hypersetup{
    colorlinks,
    linkcolor={red!50!black},
    citecolor={blue!50!black},
    urlcolor={blue!80!black}
}

\pgfplotsset{
  log x ticks with fixed point/.style={
      xticklabel={
        \pgfkeys{/pgf/fpu=true}
        \pgfmathparse{exp(\tick)}%
        \pgfmathprintnumber[fixed relative, precision=3]{\pgfmathresult}
        \pgfkeys{/pgf/fpu=false}
      }
  },
  log y ticks with fixed point/.style={
      yticklabel={
        \pgfkeys{/pgf/fpu=true}
        \pgfmathparse{exp(\tick)}%
        \pgfmathprintnumber[fixed relative, precision=3]{\pgfmathresult}
        \pgfkeys{/pgf/fpu=false}
      }
  }
}

\newcommand{\reals}{\mathbb{R}}

\newcommand{\bigo}[1]{\mathcal{O}(#1)}

\let\grad\undefined
\let\curl\undefined
\let\div\undefined
\let\nullop\undefined
\let\d\undefined
\let\tr\undefined
\DeclareMathOperator{\grad}{grad}
\DeclareMathOperator{\curl}{curl}
\DeclareMathOperator{\div}{div}
\DeclareMathOperator{\nullop}{null}
\DeclareMathOperator{\d}{d}
\DeclareMathOperator{\tr}{tr}
\DeclareMathOperator{\Exists}{\exists}
\DeclareMathOperator{\Forall}{\forall}

\newcommand{\Hgrad}{H(\grad)}
\newcommand{\Hcurl}{H(\curl)}
\newcommand{\Hdiv}{H(\div)}
\newcommand{\Ltwo}{L^2}

\newcommand{\mesh}{\mathcal{T}_h}

\newcommand{\refline}{\hat{\mathcal{I}}}

\renewcommand{\vec}[1]{\mathbf{#1}}

\DeclareMathOperator{\diag}{diag}
\DeclareMathOperator{\spn}{span}

\newcommand{\eye}{{\mathbb{I}}}

\newcommand{\be}{\vec{e}}
\newcommand{\bu}{\vec{u}}
\newcommand{\bn}{\vec{n}}

\newcommand{\bv}{\vec{v}}
\renewcommand{\bf}{\vec{f}}

\newcommand{\bx}{\vec{x}}

\newcommand{\uu}{\underline{u}}

\newcommand{\xhat}{\hat{x}}
\newcommand{\rhat}{\hat{r}}
\newcommand{\shat}{\hat{s}}
\newcommand{\Khat}{\hat{K}}
\newcommand{\Bhat}{\hat{B}}
\newcommand{\Ahat}{\hat{A}}
\newcommand{\Dhat}{\hat{D}}
\newcommand{\Shat}{\hat{S}}

\renewcommand{\P}{\mathrm{P}}
\newcommand{\DP}{\mathrm{DP}}
\newcommand{\Q}{\mathrm{Q}}
\newcommand{\DQ}{\mathrm{DQ}}

\newcommand{\NCF}{\mathrm{NCF}}
\newcommand{\NCE}{\mathrm{NCE}}

\newtheorem{remark}{Remark}

\title{Multigrid solvers for the de Rham complex with optimal complexity in polynomial degree
\thanks{Submitted to the editors November 25, 2022.
      \funding{
         PDB was supported by the University of Oxford Mathematical Institute Graduate Scholarship.
         PEF was supported by EPSRC grants EP/V001493/1 and EP/R029423/1.
      }
      This work used the ARCHER2 UK National Supercomputing Service (\url{https://www.archer2.ac.uk}).
   }
}

\author{Pablo D.\ Brubeck\thanks{
Mathematical Institute,
University of Oxford,
Oxford UK (\email{brubeckmarti@maths.ox.ac.uk})
}
\and Patrick E.\ Farrell\thanks{
Mathematical Institute,
University of Oxford,
Oxford UK (\email{patrick.farrell@maths.ox.ac.uk})
}
}

\headers{Multigrid for $p$-FEM Riesz maps}{P.~D.~Brubeck and P.~E.~Farrell}

\begin{document}

\numberwithin{equation}{section}
\maketitle

\begin{abstract}
The Riesz maps of the $L^2$ de Rham complex frequently arise as subproblems in
the construction of fast preconditioners for more complicated problems.  In
this work we present multigrid solvers for high-order finite element
discretizations of these Riesz maps with the same time and space complexity as
sum-factorized operator application, i.e., with optimal complexity in polynomial
degree in the context of Krylov methods.  The key idea of our approach is to
build new finite elements for each space in the de Rham complex with
orthogonality properties in both the $L^2$- and $H(\d)$-inner products ($\d \in
\{\grad, \curl, \div\})$ on the reference hexahedron.  The resulting sparsity
enables the fast solution of the patch problems arising in the Pavarino,
Arnold--Falk--Winther, and Hiptmair space decompositions, in the separable case.
In the non-separable case, the method can be applied to an auxiliary operator
that is sparse by construction.  With exact Cholesky factorizations of the
sparse patch problems, the application complexity is optimal but the setup
costs and storage are not. We overcome this with the finer Hiptmair space
decomposition and the use of incomplete Cholesky factorizations imposing the
sparsity pattern arising from static condensation, which applies whether static
condensation is used for the solver or not.  This yields multigrid relaxations
with time and space complexity that are both optimal in the polynomial
degree.
\end{abstract}




\section{Introduction} \label{sec:introduction}

In this paper we introduce solvers for high-order finite element discretizations
of the following boundary value problems posed on a bounded Lipschitz domain $\Omega
\subset \mathbb{R}^d$ in $d = 3$ dimensions\footnote{Our solver
strategy extends to $d \in \mathbb{N}_+$, but we describe the case $d=3$ for concreteness.}:
\begin{alignat}{5}
\beta u - \nabla \cdot \left(\alpha \nabla u\right) &= f \mbox{~in~}\Omega, \quad
&u &= 0 \mbox{~on~} \Gamma_D, \quad
&\alpha\nabla u\cdot\bn &= 0 \text{~on~} \Gamma_N; \label{eq:hgrad}\\
\beta\bu + \nabla\times \left( \alpha \nabla\times \bu\right) &= \bf \mbox{~in~} \Omega,
&\bu\times\bn &= 0 \mbox{~on~} \Gamma_D, \quad
&\alpha\nabla \times \bu \times \bn &= 0 \text{~on~} \Gamma_N; \label{eq:hcurl}\\
\beta\bu - \nabla \left(\alpha \nabla\cdot\bu\right) &= \bf \mbox{~in~} \Omega, \quad
&\bu\cdot\bn &= 0 \mbox{~on~} \Gamma_D, \quad
&\alpha\nabla\cdot\bu &= 0 \text{~on~} \Gamma_N; \label{eq:hdiv}
\end{alignat}
where $\alpha,\beta: \Omega \to \mathbb{R}_+$ are problem parameters, $\Gamma_D \subseteq \partial \Omega$, and $\Gamma_N = \partial \Omega \setminus \Gamma_D$.
For $\alpha = \beta = 1$, these equations are the so-called \emph{Riesz maps} associated with subsets of the spaces
$H(\grad, \Omega) = H^1(\Omega)$, $H(\curl, \Omega)$, and $H(\div, \Omega)$, respectively.
These function spaces are defined as:
\begin{align}
H(\grad, \Omega) &\coloneqq
\left\{
v \in L^2(\Omega): \grad{v} \in [L^2(\Omega)]^3
\right\}, \label{eq:define_hgrad} \\
H(\curl, \Omega) &\coloneqq 
\left\{
\bv\in [L^2(\Omega)]^3 : \curl{\bv} \in [L^2(\Omega)]^3
\right\}, \label{eq:define_hcurl} \\
H(\div, \Omega) &\coloneqq 
\left\{
\bv\in [L^2(\Omega)]^3 : \div{\bv} \in L^2(\Omega)
\right\}. \label{eq:define_hdiv}
\end{align}
For brevity we shall write $\Hgrad = H(\grad, \Omega)$ etc.~where there is no potential confusion.
Our problems of interest~\eqref{eq:hgrad}--\eqref{eq:hdiv}
often arise as subproblems in the construction of fast
preconditioners for more complex systems involving solution variables in \eqref{eq:define_hgrad}--\eqref{eq:define_hdiv}~\cite{hiptmair06,mardal11},
and are the canonical maps for transforming derivatives to gradients in optimization problems posed in these spaces~\cite{schwedes17}.

The spaces \eqref{eq:define_hgrad}--\eqref{eq:define_hdiv} and 
their discretizations are organized in the $L^2$ de Rham complex
\begin{equation} \label{eq:derham}
\begin{tikzcd}
  H(\grad) \arrow[r, "\grad"] \arrow[d]  & H(\curl) \arrow[r, "\curl"]
  \arrow[d] & H(\div) \arrow[r, "\div"]  \arrow[d] & L^2  \arrow[d] \\
  \Q_p \arrow[r, "\grad"] & \NCE_p \arrow[r, "\curl"] & \NCF_p \arrow[r, "\div"] & \DQ_{p-1}
\end{tikzcd},
\end{equation}
where the complex property means that the image of one map ($\grad$, $\curl$,
or $\div$) is contained in the kernel of the next, e.g., $\grad(\Hgrad) \subset
\mathrm{ker}(\curl, \Hcurl)$.  Here $\Q_p \subset \Hgrad$, $\NCE_p \subset
\Hcurl$, $\NCF_p \subset \Hdiv$, and $\DQ_{p} \subset L^2$ are piecewise
polynomial spaces of maximum polynomial degree $p$ on a mesh $\mesh$ of
tensor-product cells (hexahedra) used for the finite element discretization of
\eqref{eq:hgrad}--\eqref{eq:hdiv}. $\NCE_p$ and $\NCF_p$ are the discrete
function spaces induced by the N\'ed\'elec edge elements and face
elements~\cite{nedelec80}, respectively.

High-order finite element discretizations are well suited to exploit modern
parallel hardware architectures. They converge exponentially fast to smooth solutions and allow for
matrix-free solvers that balance the ratio of floating point operations
(flops) to energy-intensive data movement~\cite{kolev21}.

In this work we consider multigrid methods with aggressive coarsening in $p$, where all but the finest
level in the hierarchy come from the lowest-order rediscretization of the problem at $p=1$.
Standard multigrid relaxation schemes such as point-Jacobi and Gau\ss-Seidel
are not effective for high-order discretizations of these problems; these
relaxations are only effective for \eqref{eq:hgrad} at low-order\footnote{
With less aggressive $p$-coarsening, point-Jacobi/Chebyshev smoothers have
acceptable smoothing rates with simple and efficient
implementations~\cite{thompson23, phillips23}.
}, and are never
effective for \eqref{eq:hcurl} and \eqref{eq:hdiv}.  In particular, the
convergence of the multigrid scheme is not robust with respect to $\alpha$,
$\beta$, or $p$.  However, space decompositions that experimentally exhibit convergence robust
to $\alpha$, $\beta$ and $p$ are known, proposed by Pavarino~\cite{pavarino93},
Arnold, Falk \& Winther (AFW)~\cite{arnold00}, and Hiptmair \cite{hiptmair98}.
The relaxation schemes these space decompositions induce require the solution
of patchwise problems e.g., gathering all degrees of freedom (DOFs) around each vertex,
edge, or face.

Solving these patch problems becomes challenging as $p$ increases. The storage
and factorization of the patch matrices becomes prohibitively expensive, since
standard basis functions for $\Q_p, \NCE_p$, and $\NCF_p$ introduce coupling
between all interior DOFs within a cell, causing $\bigo{p^d} \times \bigo{p^d}$
dense blocks in the matrix. The Cholesky factorization of such matrices takes
$\bigo{p^{3d}}$ flops and $\bigo{p^{2d}}$ storage.  These complexities
severely limit the use of very high polynomial degrees. In this work, we will
present an alternative strategy for solving these subproblems with
$\bigo{p^{d+1}}$ flops and $\bigo{p^d}$ storage.  These complexity bounds are
optimal in the context of Krylov methods: they match the computational
complexity of applying the discretized operator via
sum-factorization~\cite{orszag80}.

\begin{figure}
\centering
\subfloat[$A_{\grad}$, standard]{
   \includegraphics[height=0.3\textwidth,valign=c]{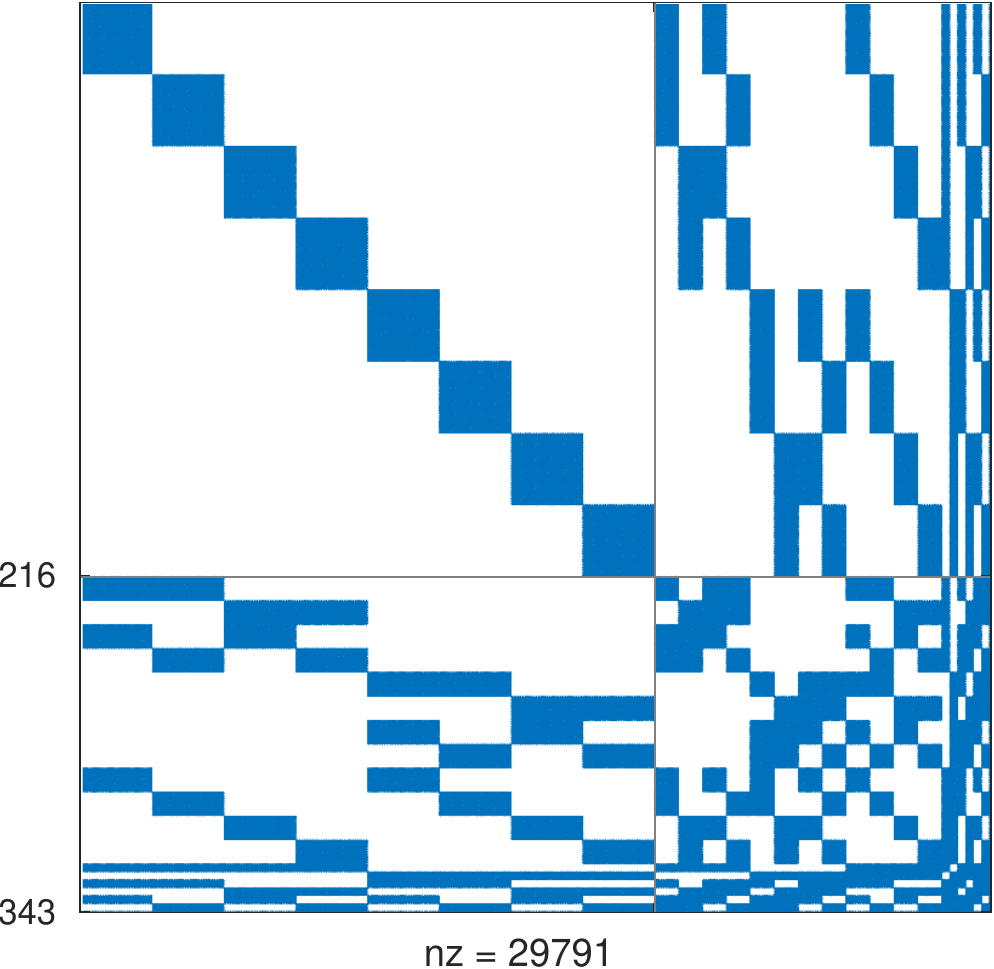}}
\hfill
\subfloat[$A_{\curl}$, standard]{
   \includegraphics[height=0.3\textwidth,valign=c]{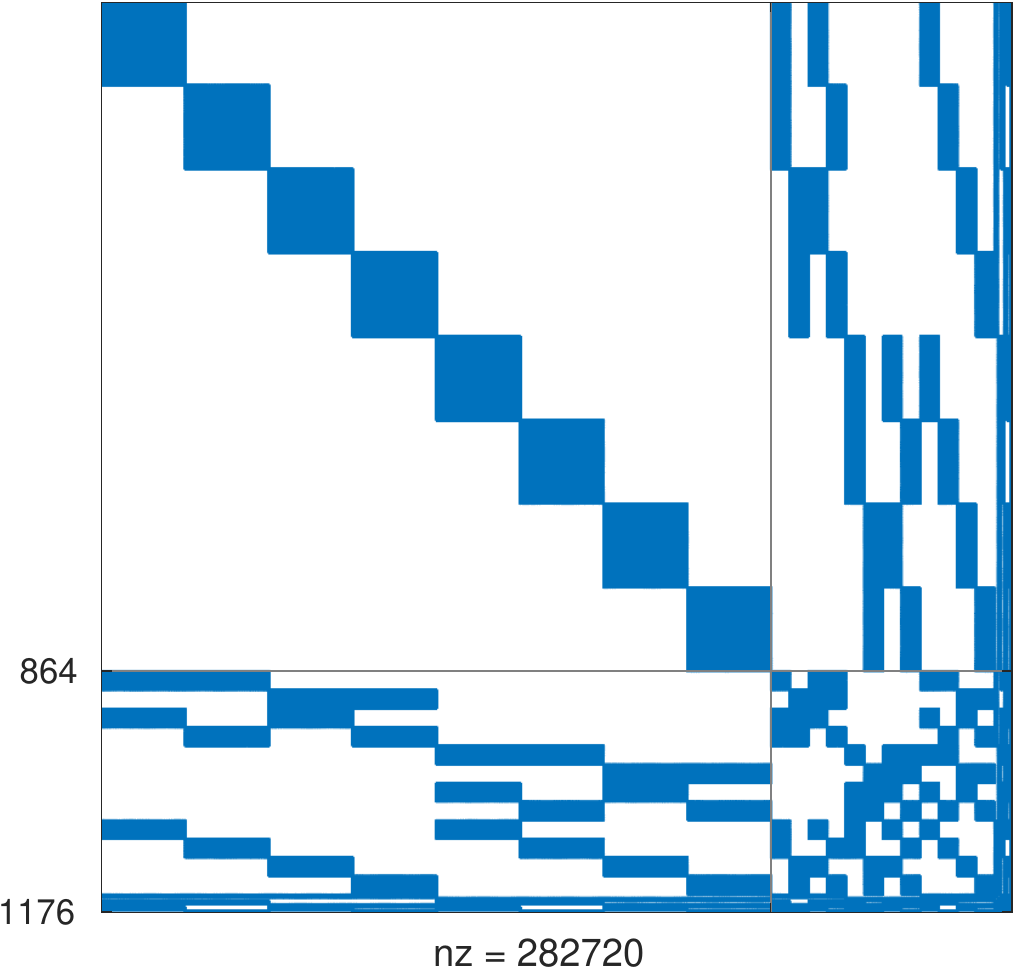}}
\hfill
\subfloat[$A_{\div}$, standard]{
   \includegraphics[height=0.3\textwidth,valign=c]{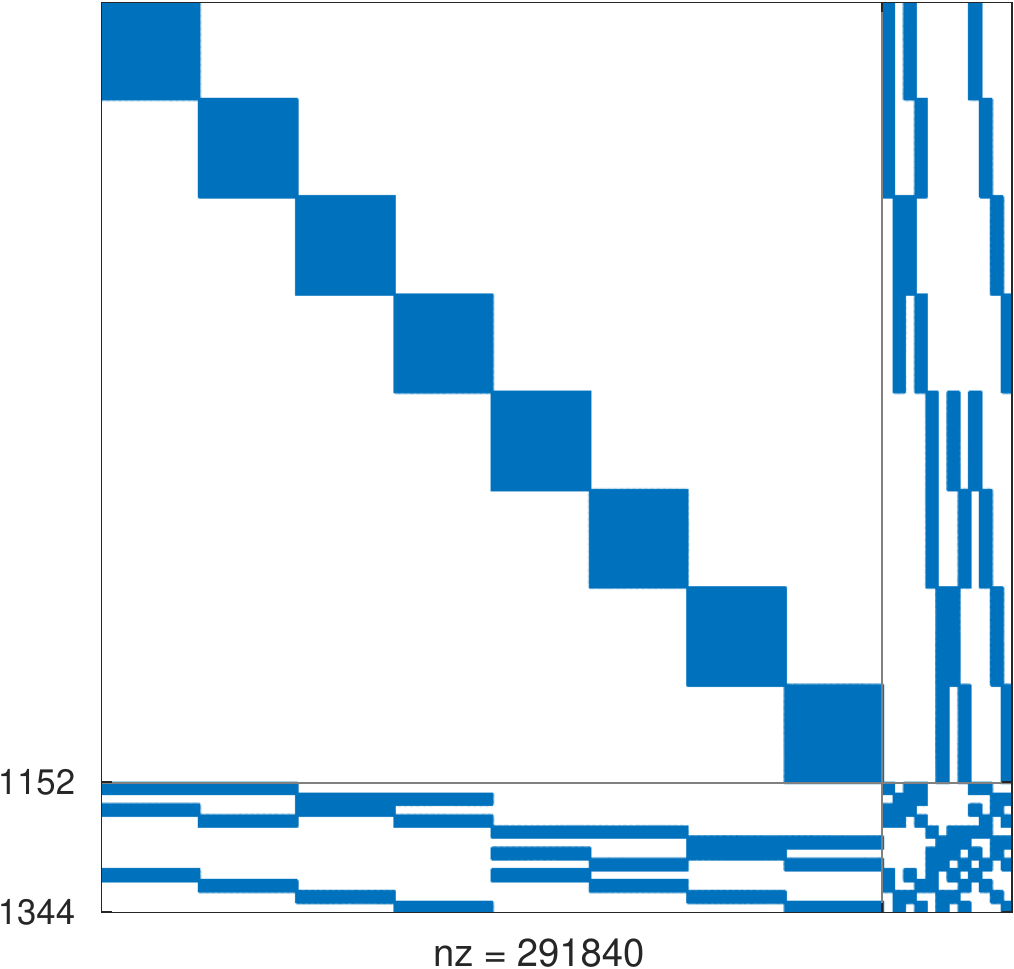}}

\subfloat[$A_{\grad}$, FDM]{
   \includegraphics[height=0.3\textwidth,valign=c]{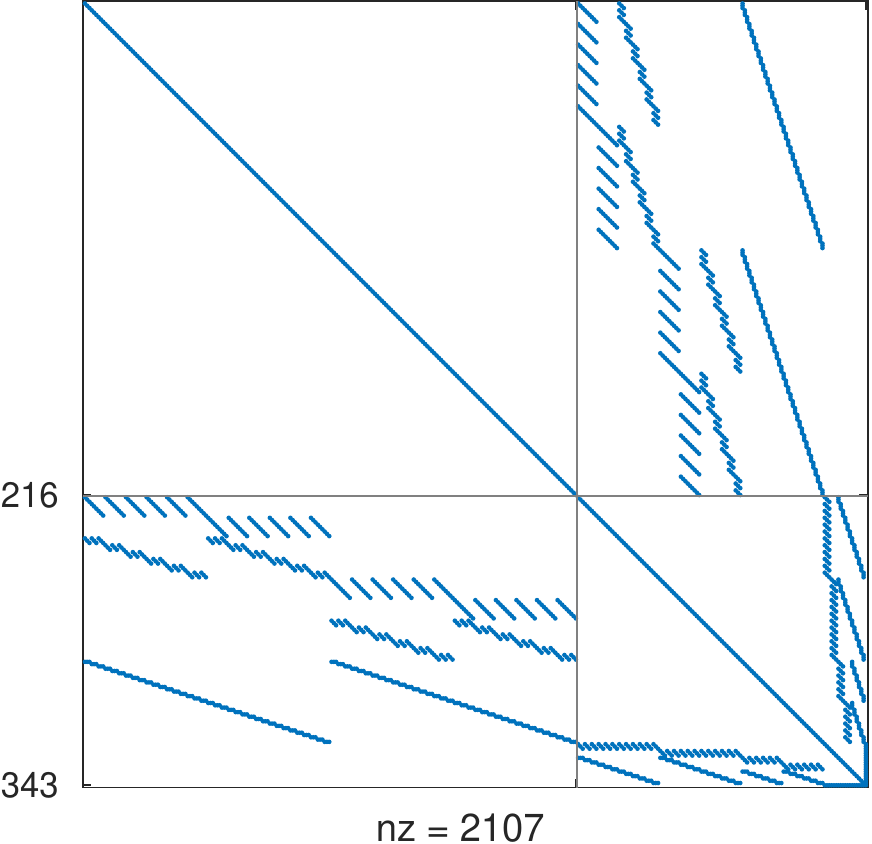}}
\hfill
\subfloat[$A_{\curl}$, FDM]{
   \includegraphics[height=0.3\textwidth,valign=c]{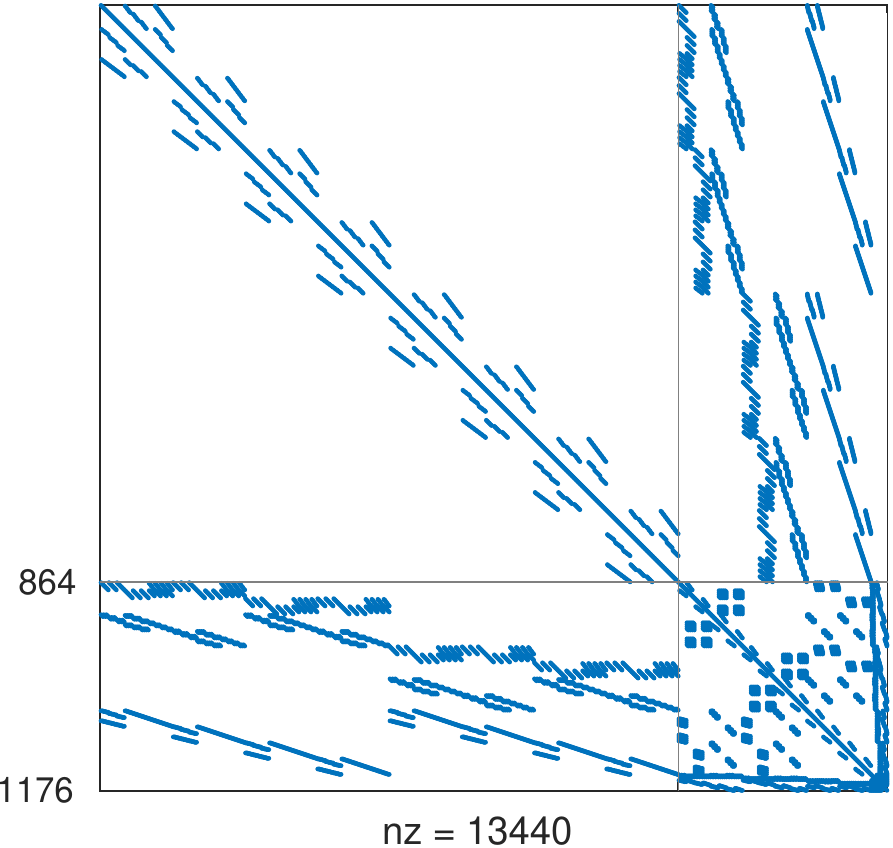}}
\hfill
\subfloat[$A_{\div}$, FDM]{
   \includegraphics[height=0.3\textwidth,valign=c]{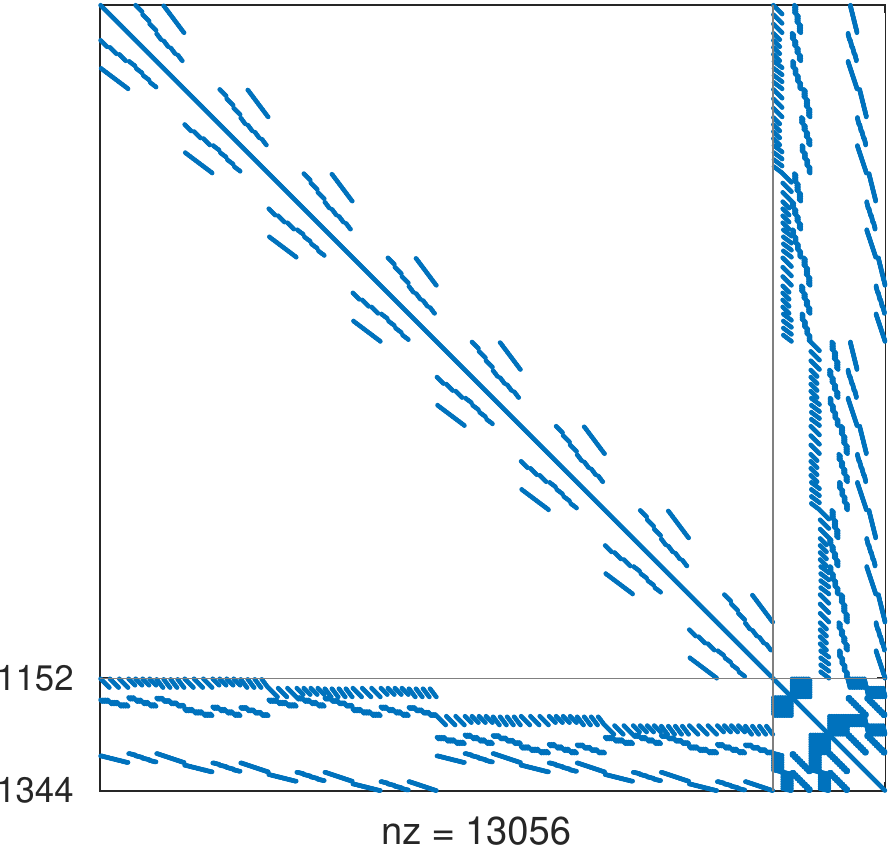}}

\caption{Sparsity patterns for the $2 \times 2 \times 2$
Pavarino--Arnold--Falk--Winther patch problem around a vertex ($p=4$), for (a, d)
\eqref{eq:hgrad} (b, e) \eqref{eq:hcurl} (c, f) \eqref{eq:hdiv}. The top row
(a)-(c) uses standard (GLL/GL) elements, while the bottom row (d)-(f) uses the
FDM variants we propose.}
\label{fig:sparsity_patterns}
\end{figure}

Our strategy relies on three main components. First, we propose new finite
elements for building
$\NCE_p$, $\NCF_p$, and $\DQ_{p}$ with useful orthogonality properties on the reference cell. These new finite elements
employ different degrees of freedom to the usual ones, and hence construct different basis functions, but discretize the same spaces. The
elements are derived from a finite element for $\Q_p$ introduced in our previous
work~\cite{brubeck22} via tensor-product construction. The new finite elements
are simple and convenient to implement; by their orthogonality properties, the
patch matrices on Cartesian cells are sparse. For example, the sparsity
patterns of a vertex-patch problem for \eqref{eq:hgrad}--\eqref{eq:hdiv} with
$p=4$ are shown in~\Cref{fig:sparsity_patterns}, for both the standard
(Gau\ss--Legendre--Lobatto/Gau\ss--Legendre, GLL/GL) finite elements and our
proposals (referred to as `FDM' elements, as they are inspired by the fast
diagonalization method~\cite{lynch64}).

The second main component is to ensure optimal fill-in in the factorization of
the patch problems. The Cholesky factorizations of the matrices shown
in~\Cref{fig:sparsity_patterns}(d-f) are sparse, even sparser than the Cholesky
factorization of a low-order discretization on a grid with the same number of
DOFs. However, this still incurs suboptimal setup and storage costs of
$\bigo{p^{2d}}$ and $\bigo{p^{d+1}}$, respectively.  We overcome this through
the choice of Hiptmair space decompositions, which require smaller patch solves
around edges and faces, and through the careful use of incomplete
factorizations of vertex patch problems.  Choosing edge patches in $\Hcurl$ and
face patches in $\Hdiv$ (along with patches for scalar and vector potential
fields, respectively) results in patch factors with fill-in of optimal space
complexity of $\bigo{p^d}$.  However, this does not address the case of
$\Hgrad$. A natural strategy is to employ incomplete Cholesky (ICC)
factorizations. The zero-fill-in ICC factorization does not work: the
factorization may fail, and even when it is computed it may not offer an
effective relaxation. Instead, we use a nested dissection ordering and impose
the sparsity pattern associated with static condensation (i.e., when the
interior DOFs are eliminated). Computational experiments indicate that this
still offers an excellent relaxation, while achieving optimality in both flops
and storage, albeit without a theoretical basis.

The third main component is the use of auxiliary operators. The patch matrices
assembled with the FDM elements are not sparse for distorted cells and/or
spatially-varying $\alpha$ or $\beta$. To overcome this, we apply our
preconditioner to an auxiliary operator which is constructed so that the patch
matrices are sparse. The auxiliary operator employed in this work is different
to that in our previous work~\cite{brubeck22} for solving \eqref{eq:hgrad}: in
this work, we construct the auxiliary operator by taking diagonal
approximations of mass matrices involved in the definition of the stiffness
matrices. The implementation of this auxiliary operator is more convenient for
$\Hcurl$ and $\Hdiv$.  The quality of this approximation depends on the mesh
distortion and the degree of any spatial variation in coefficients.
Computational experiments suggest a slow growth of the equivalence constants
with respect to $p$, but this growth is not fast enough to affect optimality of
the solver.

\begin{figure}
\centering
\input{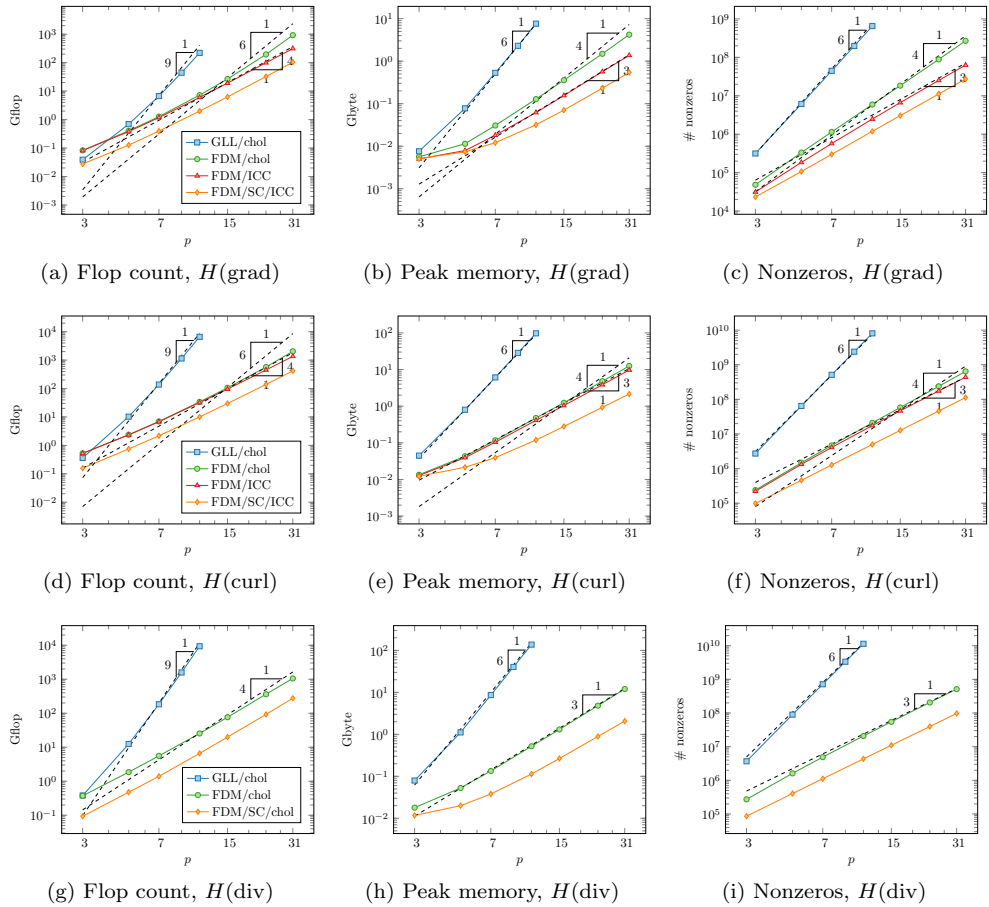}
\caption{
Flop counts, peak memory usage, and nonzeros in the sparse matrices and patch factors in 
the solution of the Riesz maps on an unstructured mesh with 27 cells.
}
\label{fig:complexities}
\end{figure}

With these components we achieve optimal complexity solvers. To illustrate
this, we show in~\Cref{fig:complexities} the number of flops and bytes required
to solve the Riesz maps \eqref{eq:hgrad}--\eqref{eq:hdiv} with conjugate
gradients (CG) and $\alpha = \beta = 1$ on an unstructured hexahedral mesh. The
setup is described in more detail in \Cref{sec:riesz_maps_experiments}.

For patch matrices assembled in the GLL basis (and hence with dense blocks) and
factorized, the time complexity is $\bigo{p^9}$, and the space complexity is
$\bigo{p^6}$, as expected.  When assembled in the FDM basis in a sparse matrix
format, the time complexity is reduced to approximately $\bigo{p^5}$ (the
empirical slope between $p = 23$ and $p = 31$ is $\bigo{p^{5.22}}$ in
\cref{fig:complexities}(a)), and the space complexity is $\bigo{p^4}$ in
\cref{fig:complexities}(c). The peak memory used, as shown in
\cref{fig:complexities}(b), is not yet scaling as $\bigo{p^4}$, indicating that
the sparse Cholesky factors are not yet the dominant term in memory.  These
complexities are further reduced when the FDM elements are combined with ICC
factorizations: the time complexity becomes $\bigo{p^4}$ and the space
complexity $\bigo{p^3}$, as desired.

A fourth, optional ingredient is the use of static condensation on the sparse
auxiliary operator, eliminating
the interior DOFs to yield patch problems posed only on the interfaces between cells.
As shown in \Cref{fig:complexities}, the time and space complexities remain
$\bigo{p^4}$ and $\bigo{p^3}$, respectively, but there are several advantages
nevertheless. First, the peak memory required is substantially reduced,
by a factor of $2.5\times$--$5.3\times$ for $p=31$ , thus enabling the use of modern HPC
hardware with limited RAM available per core. This reduction is caused by
the fact that the patch matrices have fewer DOFs. 
Second, the volume of data shared between processes is reduced, to
the same order as that required by operator application. This reduces parallel
communication. These advantages
result in a substantially faster solver.

\subsection{Related work}

The fast diagonalization method (FDM) \cite{lynch64} is a matrix factorization
that enables the direct solution of problems such as \eqref{eq:hgrad} in
optimal complexity, whenever separation of variables is applicable (i.e., only
on certain domains and for certain coefficients $\alpha, \beta$).  The FDM
breaks down the $d$-dimensional problem into a sequence of one-dimensional
generalized eigenvalue problems.  Our construction of finite elements with
orthogonality properties below in \Cref{sec:sparsity} was
inspired by the FDM.

Hientzsch~\cite{hientzsch01,hientzsch05} studied the extension of the FDM to
the $\Hcurl$ Riesz map~\eqref{eq:hcurl}.  The algorithm relies on the
elimination of one vector component.  The reduced system for the remaining
vector component can be solved directly with the FDM. With this approach, one
has to solve one generalized eigenvalue problem for every cell of the mesh,
whereas in the $\Hgrad$ case~\eqref{eq:hgrad} one solves a single generalized
eigenvalue problem on the reference cell.  Unfortunately, this strategy does not
extend well to the case with three vector components, as the nested Schur
complement is no longer a sum of three Kronecker products.

Low-order-refined (LOR) preconditioners~\cite{orszag80, deville90} are well
known to be spectrally equivalent~\cite{canuto10} to discrete high-order
operators in $\Hgrad$.
For problems in $\Hcurl$ and $\Hdiv$, LOR preconditioners have been studied in
\cite{dohrmann21, pazner22}. The auxiliary low-order problem is sparse, but
even for the $\Hgrad$ case~\eqref{eq:hgrad}, devising efficient relaxations is
challenging~\cite{pazner20}.  This is because the low-order-refined grid is
anisotropic, and pointwise smoothers on the auxiliary low-order problem become
ineffective as $p$ increases. To overcome this, Pazner~\cite{pazner20} applies
patchwise multigrid with ICC relaxations to~\eqref{eq:hgrad}, relying on a good
ordering of the DOFs. This method has not been studied for~\eqref{eq:hcurl}
or~\eqref{eq:hdiv}, to the best of our knowledge.  Instead, Pazner, Kolev \&
Dohrmann apply the AMS and ADS algebraic multigrid solvers of Hypre
\cite{falgout02,kolev09,kolev12} to the auxiliary low-order problem, which
implement the strategy proposed by Hiptmair \& Xu \cite{hiptmair07}. These LOR
approaches do not naturally combine with static condensation, since the density
of the Schur complement arising from the elimination of DOFs that are interior
to the original grid counteracts the sparsity offered by the low-order problem.

In common with the work of Sch\"oberl \& Zaglmayr~\cite{schoberl05,
zaglmayr06}, we obtain basis functions with local
complete sequence properties, i.e., on the interior of each cell the discrete
spaces form a subcomplex of the de Rham complex. Our construction of the basis
functions for $\NCE_p$, $\NCF_p$, and $\DQ_p$ from that of $\Q_p$ is the same,
but we start from a different basis for $\Q_p$; Sch\"oberl \& Zaglmayr start
with integrated Legendre polynomials on the reference interval, whereas we
employ the FDM element proposed in our previous work~\cite{brubeck22}. This
choice yields greater sparsity, because the mass matrix on the reference
interval decouples the interior DOFs.

\section{Sparsity-promoting discretization} \label{sec:sparsity}
Our goal is to construct finite elements so that the discretizations
of~\eqref{eq:hgrad}--\eqref{eq:hdiv} are sparse even at high $p$.  Specifically,
we desire that the number of nonzeros in the stiffness matrix is of the same
order as its number of rows or columns, in certain cases (Cartesian cells and
cellwise-constant coefficients).

\subsection{Exterior calculus notation and weak formulation}
To unify the discussion of \eqref{eq:hgrad}--\eqref{eq:hdiv}, we adopt the
language of the finite element exterior calculus (FEEC)~\cite{arnold18}.  We
recognize functions in $\Hgrad$, $\Hcurl$, $\Hdiv$, and $\Ltwo$ as differential
$k$-forms for $k = 0, 1, 2, 3$, respectively, writing $H(\d^k, \Omega) = H
\Lambda^k(\Omega)$, with the exterior derivative $\d^k$ corresponding to
$\d^0 = \grad$, $\d^1 = \curl$, $\d^2 = \div$, $\d^3 = \nullop$ (the zero map). We
define the spaces
\begin{equation}
V^k \coloneqq \{ v \in H \Lambda^k(\Omega) : \tr v = 0 \mbox{~on~} \Gamma_D\},
\end{equation}
where the trace operator on $\partial\Omega$ is $\tr v = v|_{\partial \Omega}$ for $0$-forms, $\tr
v = (v \times \bn)|_{\partial \Omega}$ for $1$-forms, $\tr v = (v \cdot
\bn)|_{\partial \Omega}$ for $2$-forms, and $\tr v = 0$ for
$3$-forms~\cite{arnold18}, where $\bn$ denotes the outward-facing unit normal
on $\partial\Omega$.  In FEEC notation, the discrete spaces on the bottom row of
\eqref{eq:derham} are denoted as $\Q^-_p \Lambda^k(\mesh)$.  We denote the
discrete spaces employed as $V^k_{h,p} \coloneqq \Q^-_p \Lambda^k(\mesh) \cap V^k$.

With this notation, the common weak formulation of
\eqref{eq:hgrad}--\eqref{eq:hdiv} is to find $u \in V^k$ such that
\begin{equation} \label{eq:common_weak_formulation}
a^k(v, u) \coloneqq 
(v, \beta u)_\Omega + (\d^k v, \alpha \d^k u)_\Omega = F(v) \text{ for all } v \in V^k,
\end{equation}
for $k \in \{0, 1, 2\}$, and where $(\cdot, \cdot)_\Omega$ denotes the
$L^2(\Omega)$-inner product. The discretization we consider is to find $u_{h}
\in V^k_{h,p}$ such that
\begin{equation}
a^k(v_h, u_h) = F(v_h) \text{ for all } v_{h} \in V^k_{h,p}.
\end{equation}

\subsection{Orthogonal basis for $\Hgrad$ on the interval}

In~\cite{brubeck22} the authors introduced a new finite element for $V^0_{h,p}$
on the reference interval $\refline \coloneqq [-1, 1]$, which induces basis
functions that are orthogonal in both the $\Ltwo(\refline)$- and $H(\grad, \refline)$-inner products.
The basis functions can then be extended to tensor-product cells in arbitrary
dimensions by tensor-products. We summarize the construction
of~\cite{brubeck22} here. The essential idea is to solve the one-dimensional
generalized eigenproblem: find $\{\hat{s}_i\}_{i=1}^{p-1} \subset
\P_p(\refline)$ such that
\begin{equation} \label{eq:fdmbasis}
(\shat_i, \shat_j)_{\refline} = \delta_{ij}, \quad
(\shat'_i, \shat'_j)_{\refline} = \lambda_i \delta_{ij}, \quad
\shat_i(-1)=\shat_i(1)=0, 
\quad i,j \in 1:(p-1),
\end{equation}
where $\P_p(\refline)$ is the set of polynomials of degree $p$ on $\refline$,
$a:b \coloneqq [a, b] \cap \mathbb{Z}$, and where summation is not implied.
This generalized eigenproblem \eqref{eq:fdmbasis} is solved once for a given $p$,
offline. With these functions, we define the degrees of freedom
$\{\shat^*_i\}_{i=0}^{p}$ as
\begin{equation} \label{eq:fdmdualbasis}
\shat^*_i(v) \coloneqq
\begin{cases}
v(-1), & i=0,\\
(\shat_i, v)_{\refline}, & i \in 1:(p-1),\\
v(1), & i=p.
\end{cases}
\end{equation}
The Ciarlet triple~\cite{ciarlet02} for our element is $(\refline,
\P_p(\refline), \{\shat^*_i\}_{i=0}^{p})$. The point evaluations at the vertices
guarantee $C^0(\refline)$ continuity, and hence $H(\grad, \refline)$-conformity.

The finite element induces a reference nodal basis dual to
$\{\shat^*_i\}_{i=0}^{p}$ in the usual way \cite[eq.~(3.1.2)]{brenner08}.  The basis functions associated
with $i \in 1:(p-1)$ are the generalized eigenfunctions $\shat_i$, by
construction (cf.~\eqref{eq:fdmbasis}).  It remains to determine the interface
basis functions $\shat_0, \shat_p$.  These two functions are defined via the
duality condition $\shat^*_i(\shat_j) = \delta_{ij}$, which reads
\begin{equation} \label{eq:unisolvence}
\begin{bmatrix}
\shat_0(-1) & \shat_j(-1) & \shat_p(-1)\\
(\shat_i, \shat_0)_{\refline} & (\shat_i, \shat_j)_{\refline} & (\shat_i, \shat_p)_{\refline}\\
\shat_0(1) & \shat_j(1) & \shat_p(1)
\end{bmatrix}
=
\begin{bmatrix}
1 & 0 & 0\\
0 & \delta_{ij} & 0\\
0 & 0 & 1
\end{bmatrix},
\quad i,j\in 1:(p-1).
\end{equation}
As a direct consequence, the reference mass matrix $\Bhat_{ij} = (\shat_i,
\shat_j)_{\refline}$ for $i, j \in 0:p$ will become almost diagonal, with
the only nonzero off-diagonal entries being $\Bhat_{0p} = \Bhat_{p0}$. This is
crucial for maintaining sparsity in higher dimensions, as the stiffness matrix
on Cartesian cells in higher dimensions is the Kronecker product of reference
mass and stiffness matrices.

We obtain $\{\shat_j\}$ numerically via Lagrange interpolants. We denote by
$\Shat \in \reals^{(p+1)\times (p+1)}$ the tabulation of the basis functions
onto the GLL points, i.e.\ $\Shat_{ij} = \shat_j(\hat{\xi}_i)$, such that
$\shat_j = \ell_i \Shat_{ij} $, where $\{\ell_j\}$ are the Lagrange polynomials
associated with the GLL points $\{\hat{\xi}_i\}_{i=0}^{p}$.  The matrix of coefficients
$\Shat$ is determined as follows.  Denote the interface DOFs by $\Gamma
\coloneqq \{0, p\}$, and the interior DOFs by $I \coloneqq 1:(p-1)$.  From the
first and last rows of \eqref{eq:unisolvence} we deduce that $\Shat_{\Gamma I}
= 0$ and $\Shat_{\Gamma\Gamma} = \eye$, with $\eye$ the identity matrix.  To
determine the tabulation of the interior basis functions onto
$\{\hat{\xi_i}\}_{i\in I}$, we note that \eqref{eq:fdmbasis} is equivalent to
the generalized eigenvalue problem: find $\Shat_{II}\in \reals^{(p-1)\times
(p-1)}$ and $\{\lambda_j\}_{j\in I}$ such that
\begin{equation}\label{eq:eigs}
\Shat_{II}^\top \Ahat^\mathrm{GLL}_{II} \Shat_{II} = \Lambda_{II},
\quad
\Shat_{II}^\top \Bhat^\mathrm{GLL}_{II} \Shat_{II} = \eye.
\end{equation}
Here $[\Ahat^\mathrm{GLL}]_{ij} = (\ell'_i, \ell'_j)_{\refline}$,
$[\Bhat^\mathrm{GLL}]_{ij} = (\ell_i,\ell_j)_{\refline}$ are the
stiffness and mass matrices discretized in the GLL basis, and $\Lambda_{II} =
\diag(\lambda_1,\ldots, \lambda_{p-1})$ is the diagonal matrix of eigenvalues.
We solve this problem numerically with the LAPACK routine $\texttt{dsygv}$~\cite{anderson99}, which
uses a Cholesky factorization of $\Bhat^\mathrm{GLL}_{II}$ and the QR algorithm
on a standard eigenvalue problem.

To determine
$\Shat_{I\Gamma}$, we employ the duality condition $(\shat_i, \shat_j)_{\refline} =
0$ for $i\in I, j\in \Gamma$ to obtain
\begin{equation} 
\Shat_{II}^\top (\Bhat^\mathrm{GLL}_{II} \Shat_{I\Gamma}  + \Bhat^\mathrm{GLL}_{I\Gamma} \Shat_{\Gamma\Gamma}) = 0.
\end{equation}
Using \eqref{eq:eigs} and $\Shat_{\Gamma\Gamma} = \eye$, we obtain
\begin{equation}
\Shat_{I\Gamma}  = - \Shat_{II}\Shat_{II}^\top \Bhat^\mathrm{GLL}_{I\Gamma}.
\end{equation}

With this element for $V^0_{h,p}$, discretizations of \eqref{eq:hgrad} on
Cartesian cells (axis-aligned hexahedra) are sparse, as sparse as a low-order
discretization, with a sparser Cholesky factorization. Because
of~\eqref{eq:fdmbasis}, the interior block of the stiffness matrix on such a
Cartesian cell is diagonal.

In the next subsections, we build upon the results of~\cite{brubeck22} to
construct interior-orthogonal bases for $V^k_{h,p}$, $k \in 1:3$.

\subsection{Orthogonal basis for $\Ltwo$ on the interval}

We first define a basis $\{\rhat_j\}_{j=0}^{p-1}$ for the space of
discontinuous polynomials of degree $p-1$ on the interval,
$\DP_{p-1}(\refline)$, by exploiting the fact that $\d(\P_p) = \DP_{p-1}$
(where $\d$ is the one-dimensional derivative operator).  We define the basis
for $\DP_{p-1}(\refline)$ as the derivatives of the interior basis functions
for $\P_p$ defined above, $\{\shat'_j\}_{j\in I}$, augmented with the constant
function:
\begin{equation}
\rhat_j \coloneqq 
\begin{cases}
\lambda_0^{-1/2} & j=0, \\
\lambda_j^{-1/2} \shat'_j & j=1,\ldots p-1.
\end{cases}
\end{equation}
Here $\lambda_0 \coloneqq |\refline|$, and $\lambda_j \coloneqq (\shat'_j,
\shat'_j)_{\refline}$ for $j\in I$ are required to normalize the basis. By
construction, the set $\{\shat'_j\}_{j\in I}$ is orthogonal in the
$\Ltwo(\refline)$-inner product. In addition $(\rhat_0, \rhat_j)_{\refline} =
0$ for $j\in I$, which follows from the fact that the interior basis functions
$\{\shat_j\}_{j\in I}$ vanish at the endpoints of $\refline$:
\begin{equation}
(\rhat_0, \rhat_j)_{\refline} = 
(\lambda_0\lambda_j)^{-1/2} \int_{\refline} \shat'_j \,{\d}\hat{x} =
(\lambda_0\lambda_j)^{-1/2} \left(\shat_j(1) - \shat_j(-1) \right) = 0. 
\end{equation}

\begin{figure}[tbhp]
\centering
\quad\hfill
\subfloat[$\P_4(\mathrm{FDM})$, $\hat{s}_j(\hat{x})$]{
   \includegraphics[height=0.24\textwidth,valign=c]{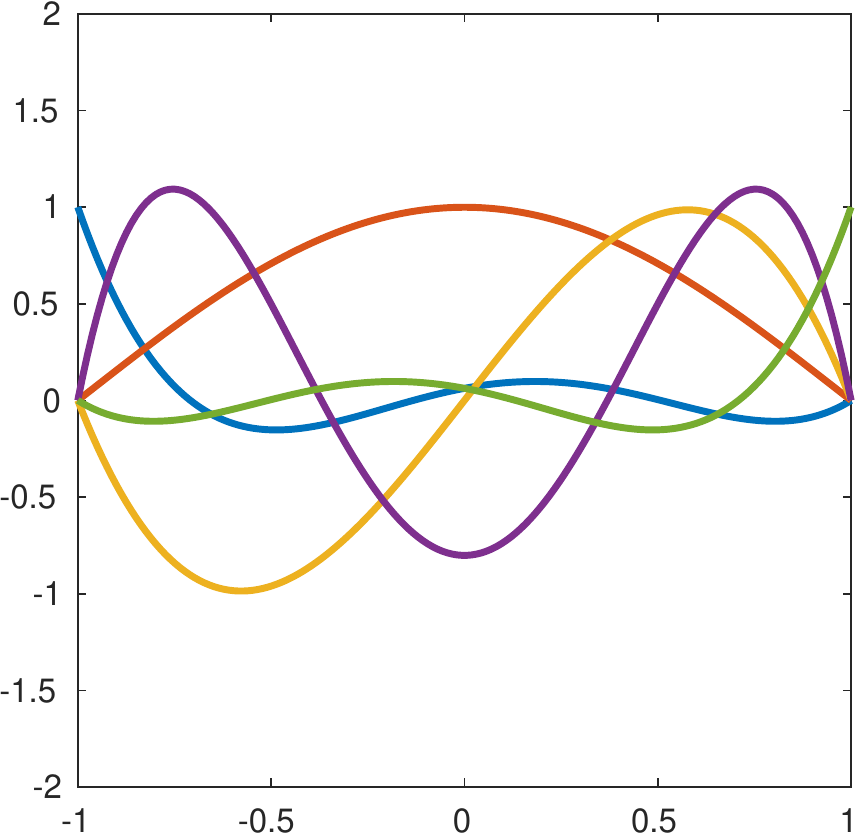}}
\hfill
\subfloat[$\DP_3(\mathrm{FDM})$, $\hat{r}_j(\hat{x})$]{
   \includegraphics[height=0.24\textwidth,valign=c]{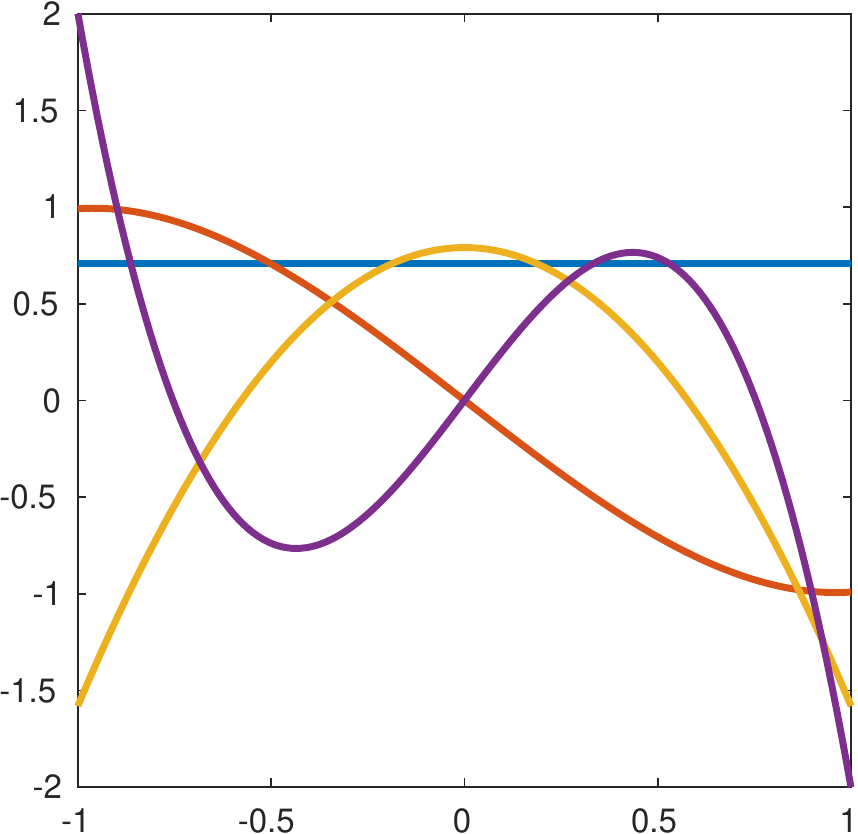}}
\hfill
\subfloat[$\hat{D}_{ij} = (\rhat_i,\shat'_j)_{\refline}$, $p=7$]{
   \includegraphics[height=0.24\textwidth,valign=c]{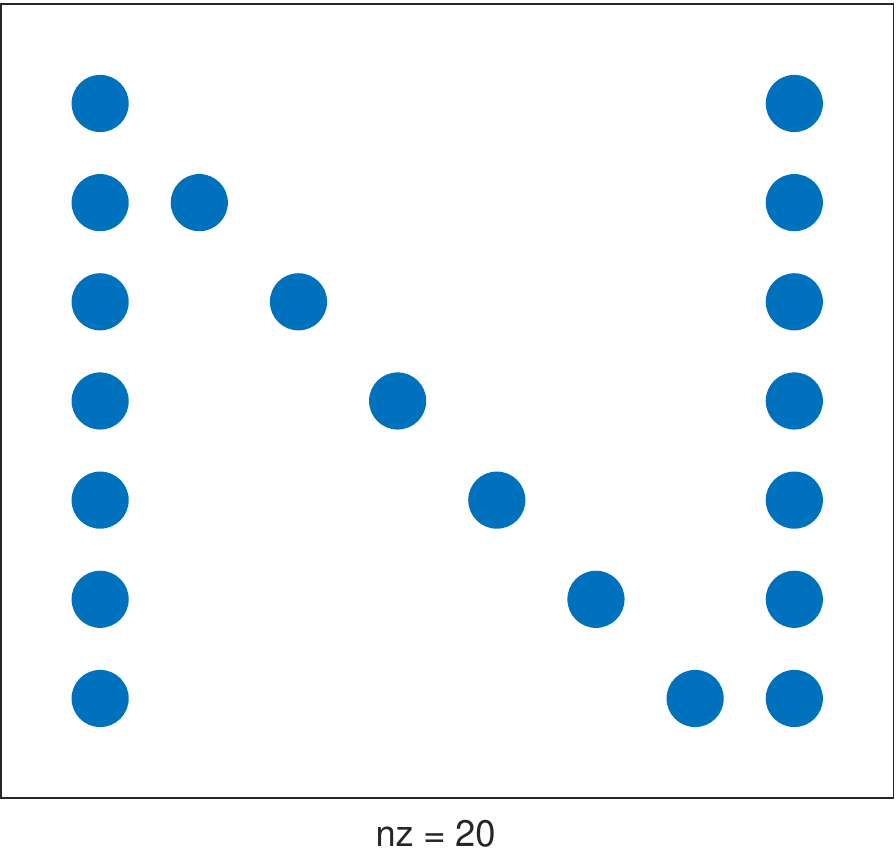}}
\hfill\quad

\caption{
Plots of the FDM basis functions for $p = 4$ on the reference interval $\refline$.}
\label{fig:fdmbasis}
\end{figure}

This dependence of the basis of $\DP_{p-1}$ on that of $\P_p$ becomes useful in
higher dimensions for enforcing interior-orthogonality in $\Hcurl$ and $\Hdiv$.
\Cref{fig:fdmbasis} shows the FDM basis functions for $\P_p$ and $\DP_{p-1}$
and the nonzero structure of the one-dimensional differentiation matrix
$\hat{D}\in \reals^{p\times (p+1)}$ that interpolates the derivatives of $\P_p$
onto $\DP_{p-1}$ in the FDM bases.

\subsection{Orthogonal bases for the de Rham complex}

With bases for $\P_{p}(\refline)$ and $\DP_{p-1}(\refline)$, we construct the
basis functions for $V^k_{h,p}(\Khat)$, $k \in 0:3$, on the reference
hexahedron $\Khat \coloneqq \refline^3$ in the usual tensor-product
fashion~\cite{nedelec80, arnold15}:
\begin{align}
V^0_{h,p}(\Khat) = \Q_{p}(\Khat) &= 
\P_p(\refline) \otimes
\P_p(\refline) \otimes
\P_p(\refline), \\
V^1_{h,p}(\Khat) = \NCE_p(\Khat) &= \begin{bmatrix}
\DP_{p-1}(\refline) \otimes \P_{p}(\refline) \otimes \P_{p}(\refline)\\ 
\P_{p}(\refline) \otimes \DP_{p-1}(\refline) \otimes \P_{p}(\refline)\\ 
\P_{p}(\refline) \otimes \P_{p}(\refline) \otimes \DP_{p-1}(\refline)
\end{bmatrix},\\
V^2_{h,p}(\Khat) = \NCF_p(\Khat) &= \begin{bmatrix}
\P_{p}(\refline) \otimes \DP_{p-1}(\refline) \otimes \DP_{p-1}(\refline)\\ 
\DP_{p-1}(\refline) \otimes \P_{p}(\refline) \otimes \DP_{p-1}(\refline)\\ 
\DP_{p-1}(\refline) \otimes \DP_{p-1}(\refline) \otimes \P_{p}(\refline)
\end{bmatrix},\\
V^3_{h,p}(\Khat) = \DQ_{p-1}(\Khat) &= 
\DP_{p-1}(\refline) \otimes
\DP_{p-1}(\refline) \otimes
\DP_{p-1}(\refline).
\end{align}

We introduce tensor-product bases for each finite element space in 
\eqref{eq:derham}. For $\Q_p(\Khat)$ we define $\{\hat{\psi}_{ijl}\}$ as
\begin{equation} 
\label{eq:hgrad-basis}
\hat{\psi}_{ijl} \coloneqq \shat_i(\xhat_1)\shat_j(\xhat_2)\shat_l(\xhat_3),
\quad (i,j,l)\in (0:p)^3.
\end{equation}
For $\NCE_p(\Khat)$ we define $\{\hat{\Psi}^{(m)}_{ijl}\}$ as
\begin{equation}
\begin{aligned}
\label{eq:hcurl-basis}
\hat{\Psi}^{(1)}_{ijl}&\coloneqq\rhat_i(\xhat_1)\shat_j(\xhat_2)\shat_l(\xhat_3) \be_1, \quad (i,j,l)\in (0:p-1)\times (0:p)\times (0:p),\\
\hat{\Psi}^{(2)}_{ijl}&\coloneqq\shat_i(\xhat_1)\rhat_j(\xhat_2)\shat_l(\xhat_3) \be_2, \quad (i,j,l)\in (0:p)\times (0:p-1)\times (0:p),\\
\hat{\Psi}^{(3)}_{ijl}&\coloneqq\shat_i(\xhat_1)\shat_j(\xhat_2)\rhat_l(\xhat_3) \be_3, \quad (i,j,l)\in (0:p)\times (0:p)\times (0:p-1).
\end{aligned}
\end{equation}
For $\NCF_p(\Khat)$ we define $\{\hat{\Phi}^{(m)}_{ijl}\}$ as
\begin{equation}
\begin{aligned}
\label{eq:hdiv-basis}
\hat{\Phi}^{(1)}_{ijl}&\coloneqq\shat_i(\xhat_1)\rhat_j(\xhat_2)\rhat_l(\xhat_3) \be_1, \quad (i,j,l)\in (0:p)\times (0:p-1)\times (0:p-1),\\
\hat{\Phi}^{(2)}_{ijl}&\coloneqq\rhat_i(\xhat_1)\shat_j(\xhat_2)\rhat_l(\xhat_3) \be_2, \quad (i,j,l)\in (0:p-1)\times (0:p)\times (0:p-1),\\
\hat{\Phi}^{(3)}_{ijl}&\coloneqq\rhat_i(\xhat_1)\rhat_j(\xhat_2)\shat_l(\xhat_3) \be_3, \quad (i,j,l)\in (0:p-1)\times (0:p-1)\times (0:p).
\end{aligned}
\end{equation}
For $\DQ_{p-1}(\Khat)$ we define $\{\hat{\phi}_{ijl}\}$ as
\begin{equation}
\label{eq:ltwo-basis}
\hat{\phi}_{ijl} \coloneqq \rhat_i(\xhat_1)\rhat_j(\xhat_2)\rhat_l(\xhat_3), 
\quad (i,j,l)\in (0:p-1)^3.
\end{equation}
By construction, the interior basis functions of these four bases are orthonormal
in the $\Ltwo(\Khat)$-inner product.
Moreover, each horizontal arrow in \eqref{eq:derham} gives rise to the following relations
between the interior basis functions, where $(i,j,l) \in (1:p-1)^3$:
\begin{equation} \label{eq:gradbasis}
\grad\hat{\psi}_{ijl} = 
\lambda_i^{1/2} \hat{\Psi}_{ijl}^{(1)} +
\lambda_j^{1/2} \hat{\Psi}_{ijl}^{(2)} +
\lambda_l^{1/2} \hat{\Psi}_{ijl}^{(3)}, 
\end{equation}
\begin{equation} \label{eq:curlbasis}
\begin{aligned}
\curl\hat{\Psi}_{ijl}^{(1)} &= \lambda_j^{1/2}\hat{\Phi}_{ijl}^{(3)} - \lambda_l^{1/2}\hat{\Phi}_{ijl}^{(2)},\\
\curl\hat{\Psi}_{ijl}^{(2)} &= \lambda_l^{1/2}\hat{\Phi}_{ijl}^{(1)} - \lambda_i^{1/2}\hat{\Phi}_{ijl}^{(3)},\\
\curl\hat{\Psi}_{ijl}^{(3)} &= \lambda_i^{1/2}\hat{\Phi}_{ijl}^{(2)} - \lambda_j^{1/2}\hat{\Phi}_{ijl}^{(1)},
\end{aligned}
\end{equation}
\begin{equation} \label{eq:divbasis}
\div\hat{\Phi}_{ijl}^{(1)} = \lambda_i^{1/2}\hat{\phi}_{ijl},\quad
\div\hat{\Phi}_{ijl}^{(2)} = \lambda_j^{1/2}\hat{\phi}_{ijl},\quad
\div\hat{\Phi}_{ijl}^{(3)} = \lambda_l^{1/2}\hat{\phi}_{ijl}.
\end{equation}
Therefore, the FDM bases form a local complete sequence
at the interior DOF level, i.e., for a fixed $(i,j,l) \in (1:p-1)^3$ we can establish
a subcomplex of the discrete de Rham complex on the reference cube
\begin{equation} \label{eq:derhamdof}
\begin{tikzcd}
  \Q_p(\Khat) \arrow[r, "\grad"] \arrow[d]  & \NCE_p(\Khat) \arrow[r, "\curl"]
  \arrow[d] & \NCF_p(\Khat) \arrow[r, "\div"]  \arrow[d] & \DQ_{p-1}(\Khat)  \arrow[d] \\
  \spn\{\hat{\psi}_{ijl}\} \arrow[r, "\grad"] 
& \spn\{\hat{\Psi}^{(m)}_{ijl}\}_{m=1}^3 \arrow[r, "\curl"] 
& \spn\{\hat{\Phi}^{(m)}_{ijl}\}_{m=1}^3 \arrow[r, "\div"] 
& \spn\{\hat{\phi}_{ijl}\}
\end{tikzcd}.
\end{equation}
Taking into account the $\Ltwo(\Khat)$-orthogonality of the interior basis
functions, \eqref{eq:derhamdof} implies that on Cartesian cells, the sparsity
pattern of the stiffness matrices $A^k$ discretizing the bilinear form for the
Riesz map $a^k(\cdot,\cdot)$ in~\eqref{eq:common_weak_formulation} connects
each interior DOF only to interior DOFs that share $(i,j,l)$. Thus the interior
block of $A^k$ has at most $d$ nonzeros per row (for $k \in \{1, 2\}$) or one
nonzero per row (for $k \in \{0, 3\}$), as depicted
in~\Cref{fig:sparsity_patterns}.

\section{Auxiliary sparse preconditioning} \label{sec:auxiliary_operator}

The orthogonality on the reference cell, and the subsequent sparsity of the
mass and stiffness matrices, will only carry over to cells that are Cartesian
and when $\alpha, \beta$ are cellwise constant. In this section we construct a
preconditioner that extends the sparsity we would have in the Cartesian case to
the case of practical interest, with distorted cells and spatially varying
coefficients. The essential idea is to build an auxiliary operator that is
sparse in the FDM basis, by construction.  To explain this, we must first
introduce some notions of finite element assembly.

\subsection{Pullbacks and finite element assembly}
The discrete spaces $V^k_{h,p}$ are defined in such way that the trace is
continuous across facets.  This is achieved through the \emph{pullback}
$\mathcal{F}^k_K : V^k(\Khat) \to V^k(K)$ that maps functions on the reference
cell $\Khat$ to functions on the physical cell $K$.  The discrete spaces are
defined in terms of the pullback,
\begin{equation}
V^k_{h,p}(\mesh) \coloneqq \left\{v_h \in V^k: 
\Forall K\in \mesh 
\Exists \hat{v}\in V^k_{h,p}(\Khat) 
\mbox{~s.t.~} v_h|_K = \mathcal{F}^k_K(\hat{v})  
\right\}.
\end{equation}

The application of the pullback to a reference function can be described as the
composition of the inverse of the coordinate mapping $F_K : \Khat \to K$ with
multiplication by a factor depending on the Jacobian of the coordinate
transformation $J_K \coloneqq \mathrm{D} F_K$.  Let $u$ be a $k$-form on $K$
mapped from $\hat{u}$ in $\Khat$. Then
\begin{equation}
u(\bx) = \mathcal{F}^k_K(\hat{u}(\hat{\bx})) = 
\begin{cases}
\hat{u}(F^{-1}_K(\bx)) & k=0,\\
J_K^{-\top} \hat{u}(F^{-1}_K(\bx)) & k=1,\\
(\det J_K)^{-1}J_K \hat{u}(F^{-1}_K(\bx)) & k=2,\\
(\det J_K)^{-1} \hat{u}(F^{-1}_K(\bx)) & k=3,
\end{cases}
\end{equation}
for $\bx \in K$ mapped from $\hat{\bx}\in \Khat$ via $\bx = F_K(\hat{\bx})$.
The pullback preserves continuity of the traces of a $k$-form across cell
facets, which is the natural continuity requirement for $\d^k$.

Another key property of the pullback is that it commutes with $\d^k$.  The
exterior derivative $\d^k u$ can be mapped from that of the reference value
$\hat{\d}{}^k \hat{u}$,
\begin{equation} \label{eq:commutativity}
\d^{k}\mathcal{F}_K^k(\hat{u})
=
\mathcal{F}_K^{k+1}(\hat{\d}{}^{k} \hat{u}),
\end{equation}
where $\hat{\d}$ is the exterior derivative with respect to the reference coordinates $\hat{\bx}$.
The pullback is incorporated in FEM by storing reference values as the DOFs in
the vector of coefficients $\uu = (\hat{u}_1, \ldots, \hat{u}_N)^\top$
representing a discrete function on a cell $K$ as
\begin{equation}
u_h|_K = \sum_{j=1}^N \hat{u}_j \mathcal{F}_K^k(\hat{\psi}^k_j),
\end{equation}
where $\hat{\psi}^k_j$ indexes the basis functions for $V^k_{h,p}(\Khat)$ defined in 
\eqref{eq:hgrad-basis}--\eqref{eq:ltwo-basis}.
The assembly of a bilinear form involves the cell matrices
\begin{equation}
[A^k_K]_{ij} = a^k(\mathcal{F}_K^k(\hat{\psi}^k_i), \mathcal{F}_K^k(\hat{\psi}^k_j)).
\end{equation}

\subsection{Construction of sparse preconditioners}
We rewrite the bilinear form $a^k(\cdot, \cdot)$ in terms of reference arguments and use the
property \eqref{eq:commutativity}, to obtain
\begin{equation} \label{eq:form-reference-values}
a^k(v_h, u_h)
= 
\left(\mathcal{F}_K^{k}(\hat{v}), \beta\, \mathcal{F}_K^{k}(\hat{u}) \right)_K +
\left(\mathcal{F}_K^{k+1}(\hat{\d}{}^k\hat{v}), \alpha\, \mathcal{F}_K^{k+1}(\hat{\d}{}^k\hat{u}) \right)_K,
\end{equation}
for $v_h, u_h \in V^k_{h,p}(K)$.
From \eqref{eq:form-reference-values} we see that the second term is an inner
product of arguments in $V^{k+1}_{h,p}(K)$.  This means that the cell matrices
can be sum-factorized in terms of the differentiation matrix $\Dhat$ acting on
reference values, and weighted mass matrices on $V^k_{h,p}$ and
$V^{k+1}_{h,p}$,
\begin{equation} \label{eq:sum-factorization}
A^k_K = M^{k}_{\beta, K} + \Dhat^\top M^{k+1}_{\alpha, K} \Dhat.
\end{equation}
Intuitively, we want each of the matrices in the sum-factorization
\eqref{eq:sum-factorization} to be sparse in order to achieve sparsity in
$A^k_K$. In higher dimensions the matrix $\Dhat$ inherits the sparsity of the
one-dimensional differentiation matrix depicted in \Cref{fig:fdmbasis}(c).  On
Cartesian cells and for cellwise constant $\alpha, \beta$, the matrices
$M_{\alpha,K}^{k+1}, M_{\beta,K}^k$ are sparse, but they are not sparse when
these conditions do not hold.

We consider $M^{k}_{\beta, K}$ first.  As $M^{k}_{\beta, K}$ discretizes a
weighted $L^2(K)$-inner product, we propose to assemble the matrix in terms of
a \textit{broken} space $\overline{V}{}^k_{h,p}$ with a fully $L^2$-orthonormal
basis. In this broken space, the mass matrix to approximate will be diagonal in
the Cartesian, constant-coefficient case. The basis for $\DP_{p-1}$ was already
$L^2$-orthogonal; we define a new basis for $\overline{\P}_p$, the broken
variant of $\P_p$, by orthogonalizing the interface basis functions with
respect to each other. The interface functions were already orthogonal to the
interior ones, which follows from the definition of the interior degrees of
freedom \eqref{eq:fdmdualbasis} and the duality condition
\eqref{eq:unisolvence}.

Let $\overline{M}{}^k_{\beta,K}$ be the weighted mass matrix in the basis for $\overline{V}{}^k_{h,p}(K)$.
Then
\begin{equation} \label{eq:mass_in_broken_basis}
M^k_{\beta, K} = \overline{G}{}^\top \overline{M}{}^k_{\beta, K} \overline{G},
\end{equation}
where $\overline{G}$ is a sparse basis transformation matrix from
$V^k_{h,p}(\hat{K})$ to $\overline{V}{}^k_{h,p}(\hat{K})$. This matrix
$\overline{G}$ is block-diagonal with one block per vector component, where
each block is a Kronecker product of identity matrices and (sparse) basis
transformation matrices from $\P_p(\refline)$ to $\overline{\P}_p(\refline)$.
The matrix $\overline{M}{}^k_{\beta, K}$ is diagonal when $K$ is Cartesian and
$\beta$ is constant, unlike $M^k_{\beta, K}$ (which is sparse, but not diagonal
in this case).  To obtain a sparse approximation to $M^k_{\beta, K}$, we simply
take the diagonal of $\overline{M}{}^k_{\beta, K}$
in~\eqref{eq:mass_in_broken_basis}.  This contrasts with taking the diagonal
directly of $M^k_{\beta, K}$, which would alter the operator even when $K$ is
Cartesian and $\beta$ is constant.

Applying the same idea to $M_{\alpha, K}^{k+1}$, we approximate the stiffness
matrix $A^k_K$ with an auxiliary matrix that is sparse on any given cell, for
any spatial variation of problem coefficients:
\begin{equation} \label{eq:sum-factorization_approx}
A^k_K \approx P^k_K \coloneqq \overline{G}{}^\top\diag(\overline{M}{}^{k}_{\beta, K})\overline{G} +
\overline{D}{}^\top \diag(\overline{M}{}^{k+1}_{\alpha, K})\overline{D},
\end{equation}
where $\overline{D} \coloneqq \overline{G} \hat{D}$.

Using this auxiliary operator ensures that the patchwise problems 
that we solve in our multigrid relaxation are sparse.
We describe these patchwise problems next. 

\section{Multigrid relaxation by subspace correction} \label{sec:relaxation}

\subsection{Notation}
We now introduce the preconditioners we use to
solve~\eqref{eq:common_weak_formulation}. We express the solvers in terms of
\emph{space decompositions}~\cite{xu92}, which we summarize briefly here. Given
a discrete space $V^k_{h,p}$, the preconditioner is induced by a particular choice of
how to write it as a sum of (smaller) function spaces:
\begin{equation} \label{eq:space_decomposition}
V^k_{h,p} = \sum_i V_i.
\end{equation}
This notation for the sum of vector spaces means that for any $v_h \in V^k_{h,p}$,
there exist $\{v_i \in V_i\}_i$ such that $v_h = \sum_i v_i$. The decomposition
is not typically unique. Given an initial guess for the solution to a
variational problem posed over $V^k_{h,p}$, the Galerkin projection of the equation
for the error is solved over each $V_i$ (additively or multiplicatively).  This
gives an approximation to the error in each subspace $V_i$, which are combined.
A cycle over each subspace constitutes one step of a subspace correction
method. For more details, see Xu~\cite{xu92}.

In order to describe the space decompositions we will use, we require some
concepts from algebraic topology.  We conceive of the mesh $\mesh$ as a
\emph{regular cell complex}~\cite{pellikka13,knepley09,logg09}. This represents
the mesh as a set of entities of different dimensions, with incidence relations
between them.  For $d=3$, the entities are vertices, edges, faces, and cells,
of dimensions $k = 0, 1, 2, 3$, respectively. The incidence relations encode the
boundary operator, relating an entity of dimension $k \ge 1$ to its bounding
sub-entities of dimension $k-1$. For example, they encode that a cell has as
its sub-entities certain faces, while a face has as its sub-entities certain
edges.  Let $E^k(\mesh)$ denote the set of entities of dimension $k$ in
$\mesh$. We also define $E^{-1}(\mesh) \coloneqq \emptyset$ for notational
convenience.

\begin{figure}
\centering
\subfloat[vertex star]{
   \includegraphics[height=0.3\textwidth,valign=c]{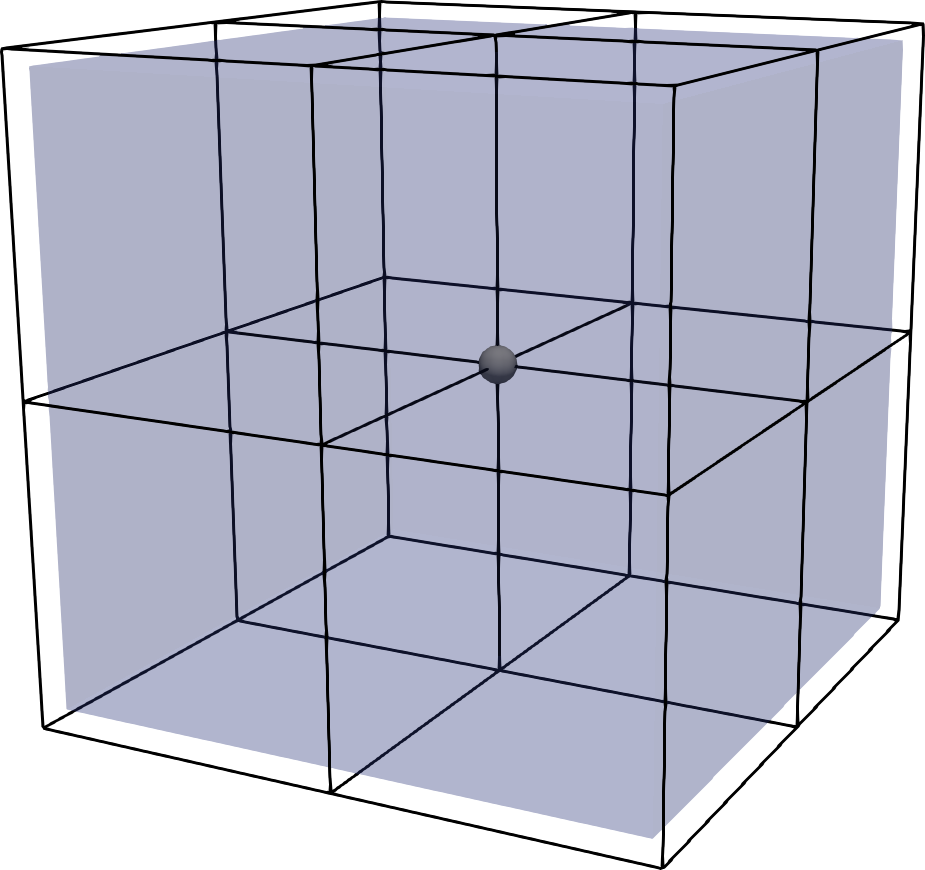}}
\hfill
\subfloat[edge star]{
   \includegraphics[height=0.3\textwidth,valign=c]{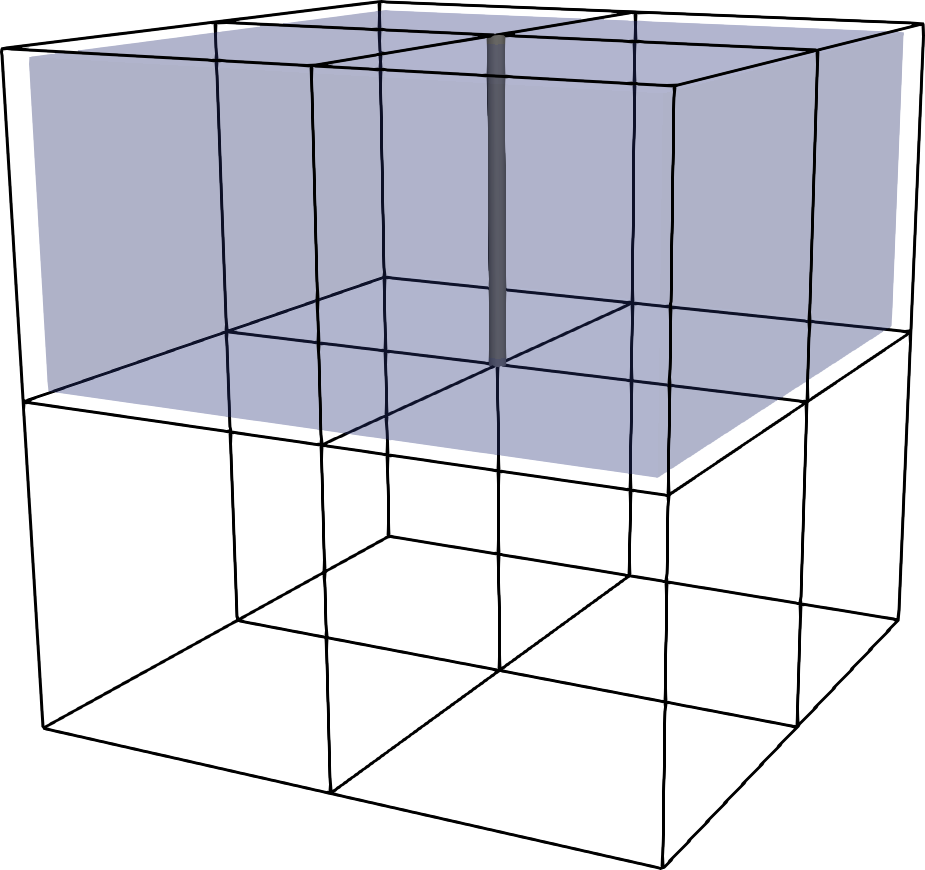}}
\hfill
\subfloat[face star]{
   \includegraphics[height=0.3\textwidth,valign=c]{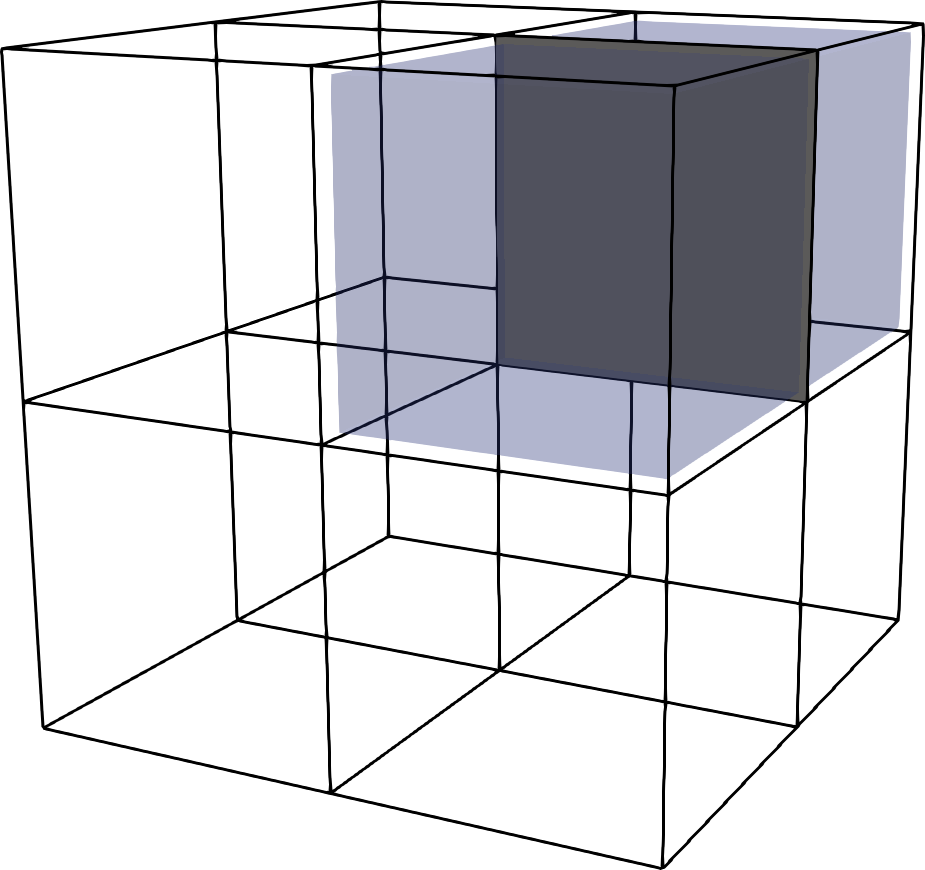}}
\caption{
The stars of (a) a vertex (b) an edge (c) a face. Solving a discrete problem
restricted to a star involves solving for all degrees of freedom contained
in the star.
}
\label{fig:stars}
\end{figure}

The \emph{star} operation on an entity $j$ of dimension $k$, denoted $\star j$,
returns the union of the interiors of all entities of dimension at least $k$ that
recursively contain $j$ as a sub-entity~\cite[\S 2]{munkres84}\cite{egorova09}.
For example, the star of a cell is simply its interior, since there are no
entities of higher dimension.  The star of an internal face $\star f$ returns the
patch of cells formed of the two cells that share $f$, excluding the
boundary of the patch.  Similarly, the star of a vertex $\star v$ returns the
union of the interiors of all edges, faces, and cells sharing $v$, as well as
the vertex itself; geometrically, this forms the patch of cells sharing $v$,
again excluding the boundary of the patch.  The stars of a vertex, edge, and
face are shown in \Cref{fig:stars}.

Given a function space $V_h$ and an entity $j$, we define
\begin{equation}
\left.V_h\right|_{\star j} \coloneqq \left\{v \in V_h : \mathrm{support}(v) \subseteq \star j\right\}.
\end{equation}
This gives a local function space around an entity upon which a variational
problem may be solved. Informally, it defines a block employed in a block
Jacobi method, taking all the DOFs in the patch of cells around
$j$, excluding those on the boundary.

\subsection{Designing space decompositions}
The framework~\eqref{eq:space_decomposition} offers a great deal of freedom in
designing solvers. Since the bilinear form $a^k(\cdot,\cdot)$ is symmetric and
coercive, powerful and general theories are available to guide the choice of
space decomposition~\cite{xu92,schoberl99,lee07}; for a summary,
see~\cite{pcpatch}. We remark on some general principles here.

First, the cost of each iteration of subspace correction will depend on the
dimensions of the subspaces $\{V_i\}$. It is therefore desirable that the space
decomposition be as fine as possible, i.e., $\mathrm{dim}(V_i)$ be as small as
possible for each $i$. If $V^k_{h,p} = \mathrm{span}\{\psi_1, \dots, \psi_J\}$,
then choosing $V_i \coloneqq \mathrm{span}\{\psi_i\}$ for $i = 1, \dots, J$
gives the finest possible space decomposition. The subspace correction methods
induced by this space decomposition are (point) Jacobi or Gau\ss--Seidel
iterations. However, as with any one-level method (without a
coarse space), the convergence of these schemes on their own is unacceptably slow.

To discuss the convergence of the scheme, we introduce some notation; we draw
this discussion from~\cite{pcpatch}. Define the operator $\mathcal{A}^k: V^k_{h,p} \to
(V^k_{h,p})^*$ associated with the bilinear form $a^k(\cdot,\cdot)$ via
\begin{equation}
\langle v, \mathcal{A}^k u\rangle \coloneqq a^k(v, u) \text{ for all } v, u \in V^k_{h,p},
\end{equation}
where $\langle \cdot, \cdot \rangle$ denotes the duality pairing.
For each subspace we denote 
the inclusion operator $\mathcal{R}_i^*: V_i \to V^k_{h,p}$ 
and its $L^2$ adjoint, 
the restriction $\mathcal{R}_i: (V^k_{h,p})^* \to V_i^*$,
and we define the restriction of $\mathcal{A}^k$ to $V_i$ by
\begin{equation}
   \langle v_i, \mathcal{A}_i u_i \rangle = \langle \mathcal{R}_i^* v_i, \mathcal{A}^k \mathcal{R}_i^* u_i \rangle \text{ for all } v_i, u_i \in V_i,
\end{equation}
i.e., $\mathcal{A}_i = \mathcal{R}_i \mathcal{A}^k \mathcal{R}_i^*$. The additive subspace correction preconditioner associated with the decomposition $\{V_i\}$ is then given by
\begin{equation}
   \mathcal{P}^{-1} = \sum_i \mathcal{R}_i^* \mathcal{A}_i^{-1} \mathcal{R}_i.
\end{equation}
Let $T = \mathcal{P}^{-1} \mathcal{A}^k$ be the preconditioned operator. Our goal is to estimate the condition number $\kappa(T)$ to bound the convergence of the conjugate gradient method. The condition number is bounded by two parameters.

The first, $N_O$, is the maximum overlap among subspaces. For each subspace
$V_i$, consider its set of overlapping subspaces $\mathrm{overlap}(V_i)
\coloneqq \{V_j: \Exists v_i \in V_i, v_j \in V_j \text{ s.t. } a(\mathcal{R}_i^* v_i, \mathcal{R}_j^*
v_j) \neq 0\}$; $N_O$ is the maximum over $i$ of $|\mathrm{overlap}(V_i)|$.
This bounds the maximal eigenvalue of $T$. It is straightforward to analyze
$N_O$ by inspection on a particular family of meshes, and $N_O$ is naturally
independent of $h$, $p$, $\alpha$, and $\beta$. The second measures the
stability of the space decomposition. Assume that there exists $c_1$ such that
\begin{equation}
\inf_{\substack{v_i \in V_i \\ \sum_i v_i = v_h}} \sum_i \|v_i\|^2_{\mathcal{A}_i} \le c_1 \|v_h\|^2_{\mathcal{A}} \text{ for all } v_h \in V^k_{h,p}.
\end{equation}
Then the minimum eigenvalue of $T$ is bounded below by $c_1^{-1}$.

In order for the convergence of our solver to be robust, we require that $c_1$
is independent of $h$, $p$, $\alpha$, and $\beta$. For example, the slow
convergence of the Jacobi iteration on the Poisson problem is because the space
decomposition is not stable in the mesh parameter $h$. This may be addressed by
incorporating a global coarse space (on a coarser mesh with $H > h$) into the
space decomposition (as in a two-level domain decomposition method, or a
multigrid method), which gives a stability constant independent of $h$ for the
Poisson equation. However, this space decomposition is still not stable in $p$
for \eqref{eq:hgrad}, nor is it stable in $h$, $\alpha$, $\beta$ for
\eqref{eq:hcurl} or \eqref{eq:hdiv}. Essentially, this is because for $k=0$,
eigenfunctions associated with small eigenvalues of the operator are smooth and
can be well-represented on a coarse grid, but this is not true for $k \in \{1,
2\}$. Building on work by Pavarino~\cite{pavarino93}, Arnold, Falk \&
Winther~\cite{arnold97,arnold00}, and Hiptmair~\cite{hiptmair98}, we now
discuss the space decompositions proposed by Sch\"oberl and
Zaglmayr~\cite{schoberl05,zaglmayr06} that are robust in $h$, $p$, $\alpha$,
and $\beta$ in numerical experiments.

\subsection{Choice of space decompositions}

The Pavarino--Arnold--Falk--Winther (PAFW) decomposition is
\begin{equation} \label{eq:pafw}
V^{k}_{h,p} = V^{k}_{h,1} + 
\sum_{v\in E^0(\mesh)} \left.V^{k\phantom{-}}_{h,p}\right|_{\star v}.
\end{equation}
This combines solving local problems on overlapping patches of cells sharing a vertex with a coarse solve at lowest order ($p=1$).
The Pavarino--Hiptmair (PH) decomposition is
\begin{equation} \label{eq:ph}
V^{k}_{h,p} = V^{k}_{h,1} + 
\sum_{i \in E^k(\mesh)} \left.V^{k\phantom{-}}_{h,p}\right|_{\star i} + 
\sum_{j \in E^{k-1}(\mesh)} \left.\d^{k-1}{} V^{k-1}_{h,p}\right|_{\star j}.
\end{equation}
These decompositions \eqref{eq:pafw} and $\eqref{eq:ph}$ coincide for $k = 0$.
For example, for $\Hcurl$ ($k=1$), the PH decomposition iterates over all edges
$i$ of the mesh, solving patch problems on the cells sharing each edge, while
$j$ iterates over all vertices.  The PH decomposition involves an auxiliary
problem on the local potential space $V^{k-1}_{h,p}|_{\star j}$:
find $\psi \in V^{k-1}_{h,p}|_{\star j}$ such that
\begin{equation}
a^k(\d^{k-1} \phi, \d^{k-1} \psi) = 
(\d^{k-1} \phi, \beta\d^{k-1} \psi)_{\star j} = L(\d^{k-1} \phi) 
\quad \text{ for all } \phi \in \left.V^{k-1}_{h,p}\right|_{\star j},
\end{equation}
since $\d^k \circ \d^{k-1} = 0$. For example, for $\Hcurl$, this becomes: for
each vertex $v$, find $\psi \in V^0_{h,p}|_{\star v}$ such that
\begin{equation*}
(\grad \phi, \beta \grad \psi)_{\star v} = L(\grad \phi) \quad \text{ for all } \phi \in \left.V^0_{h,p}\right|_{\star v},
\end{equation*}
a local scalar-valued Poisson-type problem.
\begin{remark}
In the case of $\Hdiv$ $(k=2)$, the auxiliary problem is singular, with a
kernel consisting of the curl-free functions.  One alternative is to define
an auxiliary space with this kernel removed, i.e., posed on the space
$V^{k-1}_{h,p} / \d^{k-2} V^{k-2}_{h,p}$.  Another approach is to add a
symmetric and positive-definite term on the kernel, such as $(\div\phi,
\beta \div\psi)$, as done by the Hiptmair--Xu decomposition
\cite{hiptmair07}, but with the conforming auxiliary space $[H(\grad)]^d$.
Our implementation deals with this problem in a pragmatic way by simply
adding a small multiple of the mass matrix for $V^{k-1}_{h,p}$ to the patch
matrix.
\end{remark}

\begin{remark}
The $p$-robustness of the PAFW decomposition was proven for $k = 0$ in the
tensor-product case by Pavarino~\cite{pavarino93} and in the simplicial case by
Sch\"oberl et al.~\cite{schoberl08}.  A similar decomposition (with geometric
multigrid coarsening, not coarsening in polynomial degree) was proposed by
Arnold, Falk, \& Winther for $k = 1$ and $k = 2$~\cite{arnold97,arnold00}, and
proven to be robust to mesh size $h$ and variations in constant $\alpha$ and
$\beta$. To the best of our knowledge the $p$-robustness of \eqref{eq:pafw} and
\eqref{eq:ph} have not been proven for $k = 1$ or $k = 2$, although numerical
experiments indicate that they are $p$-robust for $p \le 31$. 
\end{remark}

To simplify notation, we drop the superscript $k$ from the stiffness matrices
$A^k$.  The algebraic realization of the PAFW relaxation (combined additively)
reads
\begin{equation} \label{eq:asm-afw}
P^{-1}_\mathrm{PAFW} = R_0^\top A_0^{-1} R_0 +
\sum_{v \in E^0(\mathcal{T}_h)} R_v^\top A_v^{-1} R_v.
\end{equation}
Here $R_0$ is the restriction matrix from $(V^{k}_{h,p})^*$ to
$(V^{k}_{h,1})^*$, $A_0$ is the stiffness matrix for the original bilinear form
$a^k(\cdot, \cdot)$ rediscretized with the lowest-order element
$(p=1)$\footnote{We apologize for this notation; the use of a subscript $_0$ to
indicate the coarse grid is widely used in the domain decomposition literature.
In our case the coarse grid is formed with $p=1$, not $p=0$.}, $R_v$ are
Boolean restriction matrices onto the DOFs of each vertex-star patch $\star v$,
and $A_v = R_v A R_v^\top$ are sub-matrices of $A$ corresponding to the rows
and columns of DOFs of the patch. Furthermore, another subspace correction
method such as geometric or algebraic multigrid, may be employed as the $p$-coarse 
solver to approximate $A_0^{-1}$.

Similarly, the additive PH relaxation is implemented as
\begin{equation} \label{eq:asm-hiptmair}
P^{-1}_\mathrm{PH} = R_0^\top A_0^{-1} R_0 + 
\sum_{i \in E^k(\mathcal{T}_h)} R_i^\top A_i^{-1} R_i +
D \left(\sum_{j \in E^{k-1}(\mathcal{T}_h)} R_j^\top B_j^{-1} R_j\right) D^\top.
\end{equation}
where $R_0, A_0$ have the same meaning as in \eqref{eq:asm-afw},
$D$ is the matrix tabulating the differential operator $\d^{k-1} : V^{k-1}_{h,p} \to V^k_{h,p}$,
$R_i, R_j$ are Boolean restriction matrices onto star patches on entities of dimension 
$k$ and $k-1$, respectively,
$A_i = R_i A R_i^\top$ are patch matrices,  
and $B_j$ are patch matrices extracted from $B$, which is obtained as
the discretization of $a^k(\d^{k-1}\phi, \d^{k-1}\psi)$.
With the FDM basis it is feasible to compute and store the $D$ matrix, as it is
sparse and applies to reference values.

\subsection{Statically-condensed space decompositions}

As seen from \Cref{fig:sparsity_patterns}, the minimal coupling between
interior DOFs that arises from the orthogonality of the FDM elements invites
the use of static condensation.  Static condensation yields a finer space
decomposition with smaller subspaces by eliminating the interior DOFs.  The
overlapping subspaces only involve interface DOFs.

The statically-condensed Pavarino--Arnold--Falk--Winther (SC-PAFW) decomposition is
\begin{equation} \label{eq:sc-pafw}
V^{k}_{h,p} = V^{k}_{h,1} +
\left.V^{k\phantom{-}}_{h,p}\right|_\mathcal{I} + 
\sum_{v\in E^0(\mesh)} \left.\tilde{V}^{k\phantom{-}}_{h,p}\right|_{\star v},
\end{equation}
Here $\mathcal{I} \coloneqq E^d(\mesh)$ is the set of cell interiors of $\mesh$,
and
$V^k_{h,p}|_\mathcal{I} \coloneqq \sum_{c\in\mathcal{I}} V^k_{h,p}|_c$ is the space of
discrete functions supported on cell interiors.
The cell-interior problems do not overlap with each other and can be solved independently.
We denote by $\tilde{V}^{k}_{h,p} \coloneqq (V^k_{h,p}|_\mathcal{I})^\perp$ 
the space of discrete harmonic functions of $V^{k}_{h,p}$,
defined as the $a^k(\cdot, \cdot)$-orthogonal complement of 
$V^k_{h,p}|_\mathcal{I}$,
\begin{equation}
\tilde{V}^{k}_{h,p} = \left\{ \tilde{v}\in V^{k}_{h,p} : 
\,a^k(w, \tilde{v})=0
\Forall w \in \left.V^k_{h,p}\right|_{\mathcal{I}} 
\right\}.
\end{equation}
By definition $V^{k}_{h,p} = V^{k}_{h,p}|_{\mathcal{I}} \oplus \tilde{V}^{k}_{h,p}$;
this orthogonality leads to reduced overlap of the star patches, compared to
the non-statically-condensed case.

The statically-condensed Pavarino--Hiptmair (SC-PH) decomposition is
\begin{equation} \label{eq:sc-ph}
V^{k}_{h,p} = V^{k}_{h,1} + 
\left.V^{k\phantom{-}}_{h,p}\right|_\mathcal{I} + 
\sum_{i\in E^k(\mesh)} \left.\tilde{V}^{k\phantom{-}}_{h,p}\right|_{\star i} + 
\sum_{j\in E^{k-1}(\mesh)} \left.\d^{k-1} \tilde{V}^{k-1}_{h,p}\right|_{\star j}.
\end{equation}
These decompositions again coincide for $k = 0$.

Denote by $I, \Gamma$ the sets of interior and interface (vertex, edge, and
face) DOFs, respectively. Reordering the DOFs of $A$ yields a $2 \times 2$
block matrix, with inverse obtained from its block $LDL^\top$ decomposition
\begin{equation}
\label{eq:LDU}
\begin{bmatrix}
A_{II} & A_{I\Gamma} \\
A_{\Gamma I} & A_{\Gamma \Gamma}
\end{bmatrix}^{-1}
=
\begin{bmatrix}
\eye & -A_{II}^{-1}A_{I \Gamma}\\
 0 & \eye
\end{bmatrix}
\begin{bmatrix}
A^{-1}_{II} & 0\\
0 & S^{-1}
\end{bmatrix}
\begin{bmatrix}
\eye & 0\\
-A_{\Gamma I}A_{II}^{-1} & \eye
\end{bmatrix},
\end{equation}
where $S$ denotes the interface Schur complement
\begin{equation}
S = A_{\Gamma\Gamma} - A_{\Gamma I} A_{II}^{-1} A_{I\Gamma}.
\end{equation}
When solvers for $A_{II}$ and $S$ are available, the application of $A^{-1}$
times a residual vector can be performed with a single application of $S^{-1}$
and only two applications of $A_{II}^{-1}$. This is because the second and
third instances of $A_{II}^{-1}$ in the RHS of \eqref{eq:LDU} act on the same
interior DOFs of the incoming residual vector.

Another way to write \eqref{eq:LDU} gives rise to the additive interpretation
of the harmonic subspace correction step
\begin{equation} \label{eq:schur}
A^{-1} = R_I^\top A_{II}^{-1} R_I +
R_\Gamma^{\top} S^{-1} R_\Gamma,
\end{equation}
where $R_I$ is a Boolean restriction onto the interior DOFs, and
\begin{equation}
R_\Gamma = 
\begin{bmatrix}
- A_{\Gamma I}A_{II}^{-1} & \eye
\end{bmatrix}
\end{equation}
is the ideal restriction operator onto the discrete harmonic subspace.
$R_\Gamma^\top$ maps vectors of coefficients in $\tilde{V}^k_{h,p}$
to $V^k_{h,p}$. The orthogonality between 
$\tilde{V}^k_{h,p}$  and $V^k_{h,p}|_\mathcal{I}$
is reflected by the identity $R_\Gamma A R_I^\top = 0$. 

\begin{remark}
For Cartesian cells, in $H(\grad)$ $(k=0)$, the interior DOFs are fully
decoupled and $A_{II}$ is diagonal, while for $k=1,2$, the cell-interior problems
only couple at most $d$ DOFs, as shown in \eqref{eq:derhamdof}.  There
exists a reordering of the interior DOFs for which $A_{II}$ becomes block
diagonal with diagonal blocks of dimension at most $d\times d$, implying
that $A_{II}$ shares its sparsity pattern with its inverse.  Therefore
$A_{II}$ coincides with its zero-fill-in incomplete Cholesky factorization.
Hence we may assemble and store the Schur complement $S$, even for very high
$p$.  This also holds for the auxiliary sparse operator in
\Cref{sec:auxiliary_operator} on general cells, as it inherits the sparsity
pattern of the Cartesian case by construction.
\end{remark}

\begin{remark}
We choose to use a Krylov method on $A$, as opposed to implementing one on the
condensed system involving $S$.  The action of the true $A_{II}^{-1}$
involves the iterative solution of local problems on cell-interiors,
inducing $\bigo{p^{d+1}}$ computational cost in the application of the true $S$.
Although the sum-factorized application of the true $A_{\Gamma\Gamma}, A_{I \Gamma},
A_{\Gamma I}$ only involves $\bigo{p^d}$ flops, and the conditioning of the
preconditioned operator is generally better, this results in a longer
runtime when compared to using a Krylov method on $A$ with the
statically-condensed preconditioner built from the sparse auxiliary operator, especially for the case $k=1$.
\end{remark}

The algebraic realization of the SC-PAFW relaxation \eqref{eq:sc-pafw} approximates 
the Schur complement in \eqref{eq:schur} with
\begin{equation} \label{eq:sc-asm-afw}
S^{-1} \approx S^{-1}_{\mathrm{PAFW}} \coloneqq 
\tilde{R}_0^\top A_0^{-1} \tilde{R}_0 + 
   \sum_{v \in E^0(\mathcal{T}_h)} {\tilde{R}_v}^\top S_{v}^{-1} \tilde{R}_v,
\end{equation}
where $\tilde{R}_0$ is the restriction matrix from $(\tilde{V}^{k}_{h,p})^*$ to
$({V}^{k}_{h,1})^*$, and $\tilde{R}_v$ are the Boolean restriction matrices
onto the interface DOFs of the vertex-star patch $\star v$, and $S_{v} =
\tilde{R}_v S \tilde{R}_v^\top$.

Similarly, the SC-PH relaxation \eqref{eq:sc-ph} is implemented as
\begin{equation} \label{eq:sc-asm-hiptmair}
S^{-1}_{\mathrm{PH}} \coloneqq 
\tilde{R}_0^\top A_0^{-1} \tilde{R}_0 + 
   \sum_{i \in E^k(\mathcal{T}_h)} {\tilde{R}_i}^\top S_i^{-1} \tilde{R}_i +
\tilde{D}\left(
   \sum_{j \in E^{k-1}(\mathcal{T}_h)} {\tilde{R}_j}^\top \tilde{B}_j^{-1} \tilde{R}_j
\right) \tilde{D}^\top,
\end{equation}
where $\tilde{R}_0, A_0$ have the same meaning as in \eqref{eq:sc-asm-afw},
$\tilde{D} = D_{\Gamma\Gamma}$ is the matrix tabulating the differential 
operator $\d^{k-1} : \tilde{V}^{k-1}_{h,p} \to \tilde{V}^k_{h,p}$,
$\tilde{R}_i, \tilde{R}_j$ are Boolean restriction matrices onto the interface DOFs 
of star patches on entities of dimension $k$ and $k-1$, respectively,
$S_i = \tilde{R}_i S \tilde{R}_i^\top$ are patch matrices,  
and $\tilde{B}_j$ are patch matrices extracted from the interface Schur complement of 
the discretization of $a^k(\d^{k-1}\phi, \d^{k-1}\psi)$.

\subsection{Achieving optimal fill-in}

\begin{figure}[tbhp]
\centering
\subfloat[$\mathrm{chol}(A^0)$, $\Q_4(\star v)$]{
   \includegraphics[height=0.3\textwidth,valign=c]{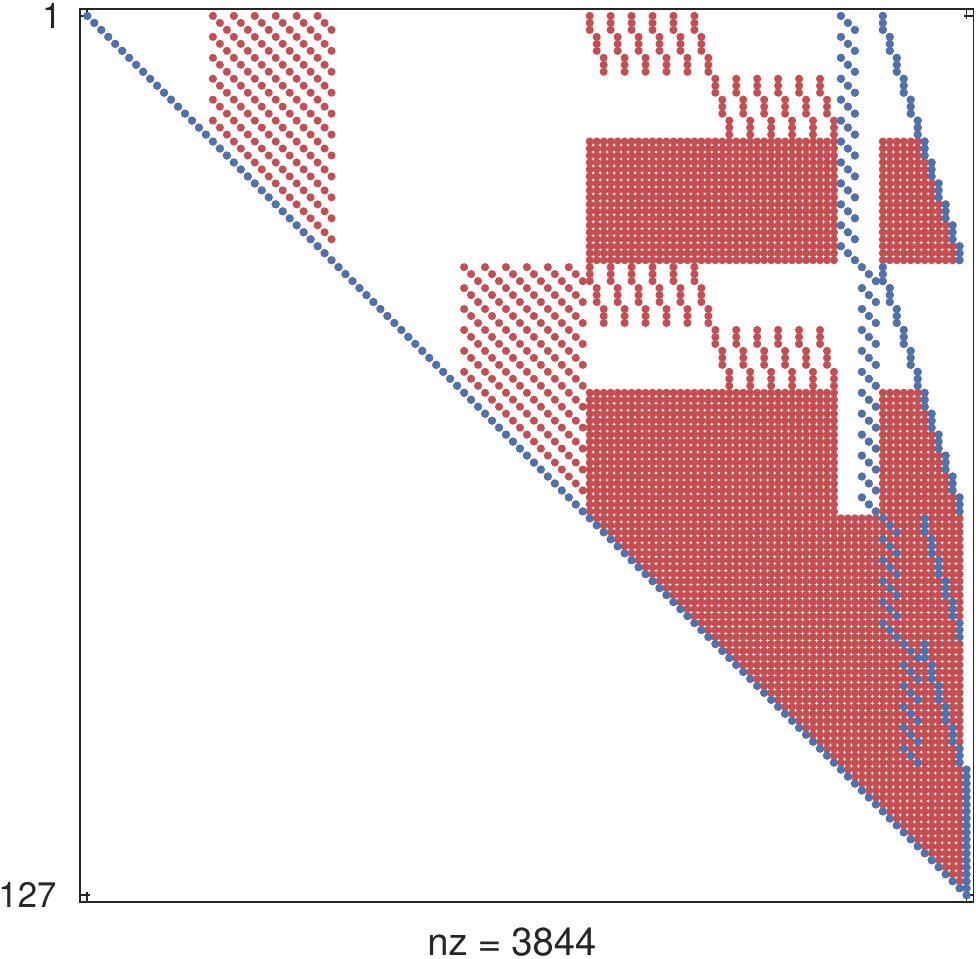}}
\hfill
\subfloat[$\mathrm{chol}(A^1)$, $\NCE_4(\star v)$]{
   \includegraphics[height=0.3\textwidth,valign=c]{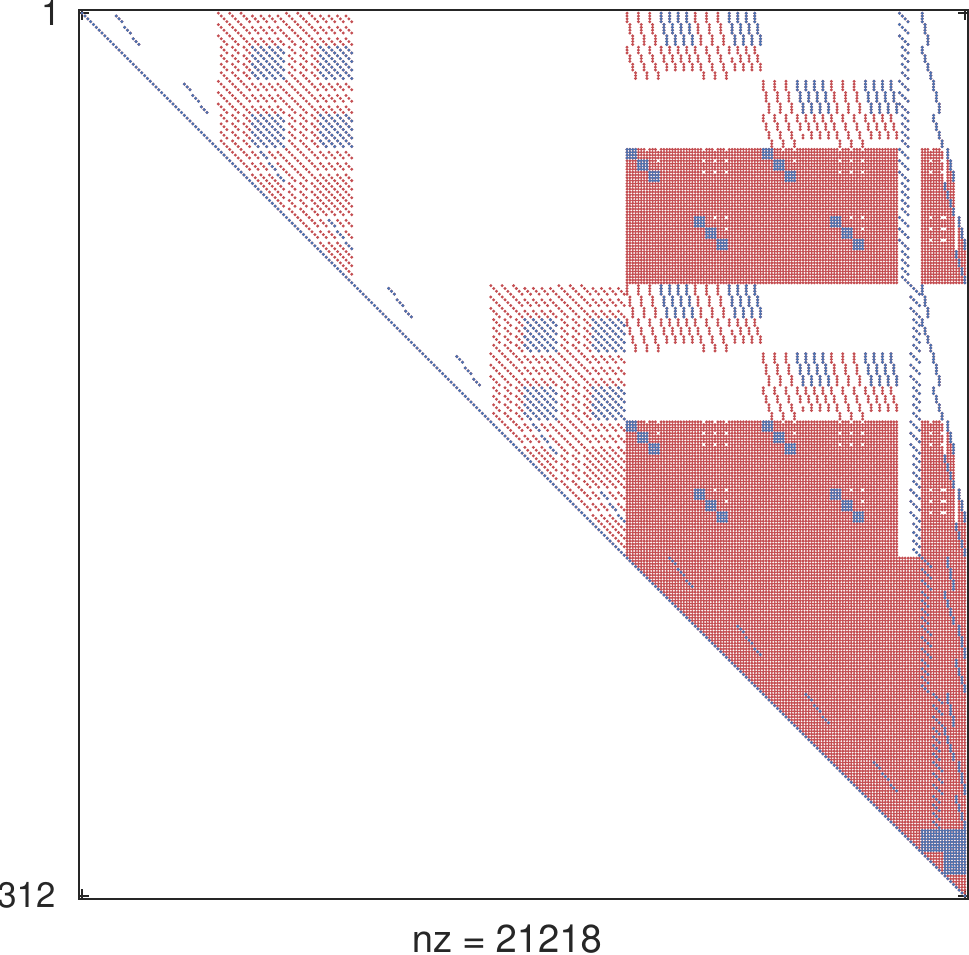}}
\hfill
\subfloat[$\mathrm{chol}(A^2$), $\NCF_4(\star v)$]{
   \includegraphics[height=0.3\textwidth,valign=c]{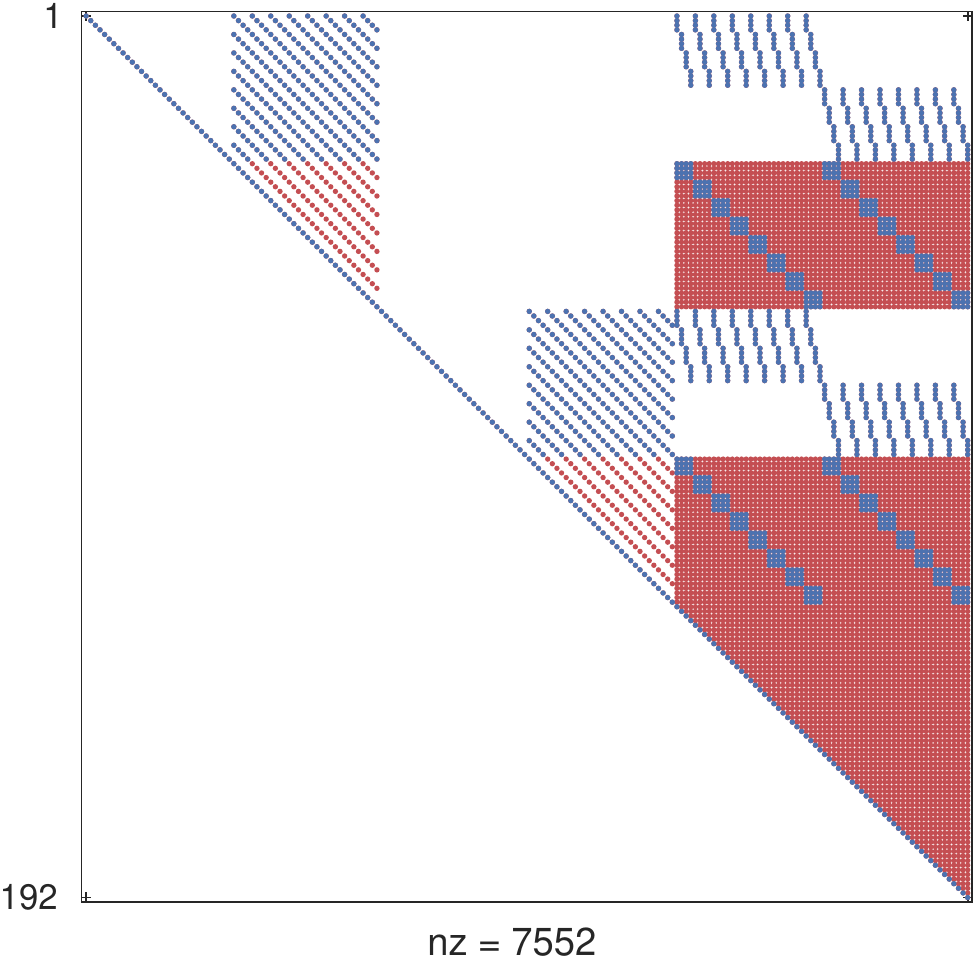}}

\subfloat[$\mathrm{ICC}(A^0)$, $\Q_4(\star v)$]{
   \includegraphics[height=0.3\textwidth,valign=c]{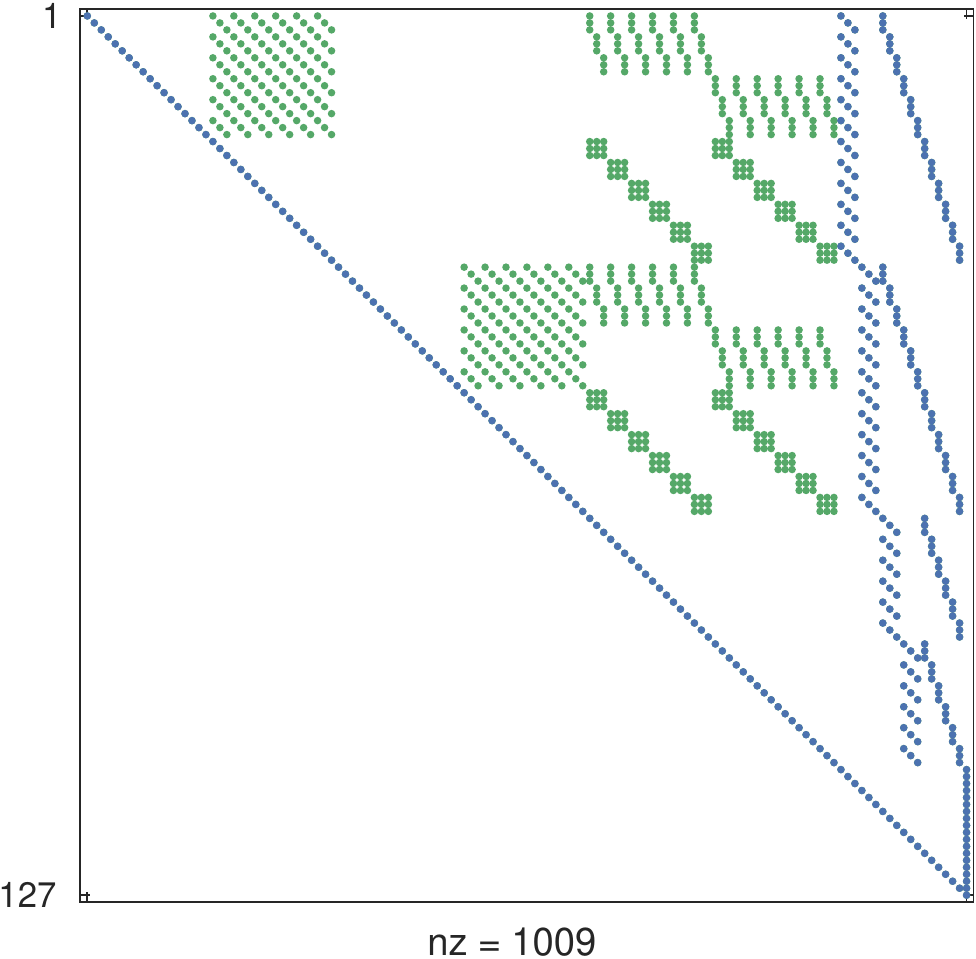}}
\hfill
\subfloat[$\mathrm{chol}(A^1)$, $\NCE_4(\star e)$]{
   \includegraphics[height=0.3\textwidth,valign=c]{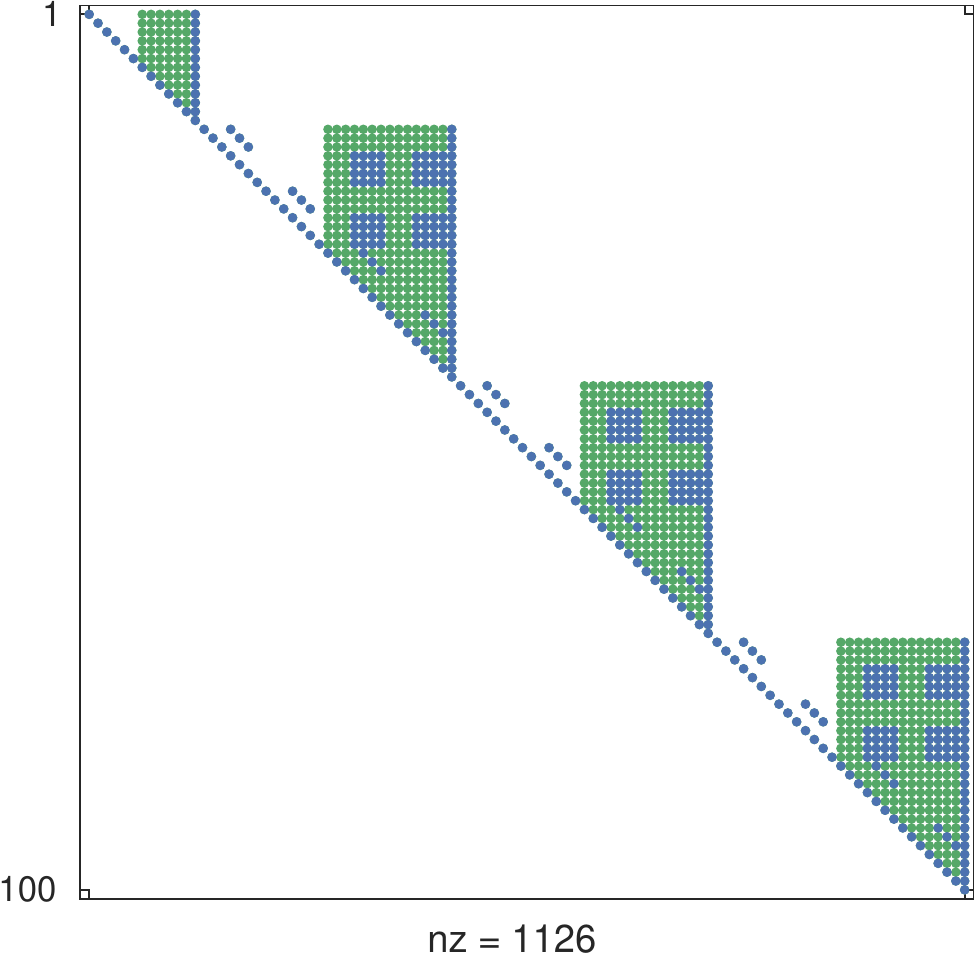}}
\hfill
\subfloat[$\mathrm{chol}(A^2)$, $\NCF_4(\star f)$]{
   \includegraphics[height=0.3\textwidth,valign=c]{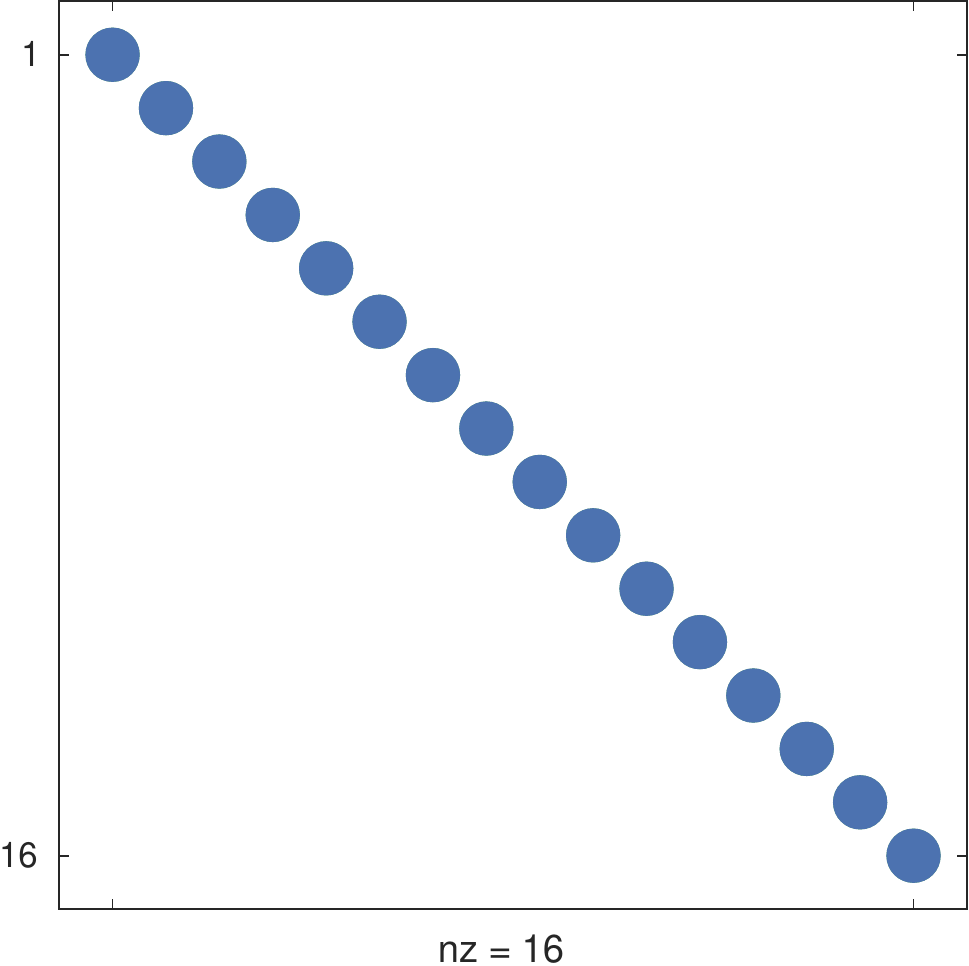}}

\caption{
The interface-interface block of the patch problems arising in the PAFW (a-c)
and the PH (d-f) space decompositions on a regular Cartesian mesh with FDM
elements. Here $\star v$, $\star e$, $\star f$ correspond to the vertex, edge, and
face stars depicted in \Cref{fig:stars}.
The sparsity pattern of the problem posed in the FDM basis is
coloured in blue. The additional nonzeros required by the Cholesky
factorization are colored in red or green, depending on whether they require
suboptimal storage (as in (a-c)) or optimal storage (as in (e-f)). In (d) we
plot the sparsity pattern of the interface Schur complement arising in
static condensation, which offers a sparsity pattern intermediate between
the zero fill-in and the full Cholesky factorization.
}
\label{fig:spy_feec}
\end{figure}

To achieve optimal complexity of our solver, we require that the factorizations
of the patch matrices arising in the space decomposition be optimal in storage.
The number of nonzeros in the factorizations also relates to the number of
flops required to compute them.  For the PAFW space decomposition
\eqref{eq:pafw}, in \Cref{fig:spy_feec}(a-c) we observe fill-in of $\bigo{p^4}$
nonzeros in the Cholesky factorization of the vertex patch problems in
$\Hgrad$, $\Hcurl$, and $\Hdiv$, even with a nested dissection ordering.
Moreover, computing this factorization incurs $\bigo{p^6}$ flops. These costs
compare unfavorably with the $\bigo{p^3}$ storage and $\bigo{p^4}$ flops
required by sum-factorized operator application.

To overcome this, for $\Hgrad$ we employ an incomplete Cholesky factorization.
Incomplete Cholesky factorization with zero fill-in (ICC(0), depicted in blue
in \Cref{fig:spy_feec}(a) and (d)) does not yield an effective relaxation
method, even when it can be computed. In contrast, the incomplete Cholesky
factorization on the statically-condensed sparsity pattern does yield an
effective relaxation. This factorization has $\bigo{p^3}$ fill-in ($\bigo{p}$
nonzeros on $\bigo{p^2}$ rows). It appears that this sparsity pattern (depicted
in green in \Cref{fig:spy_feec}(d)) offers a suitable intermediate between the
zero-fill-in pattern and the full Cholesky factorization (depicted in red in
\Cref{fig:spy_feec}(a)).

For $\Hcurl$ and $\Hdiv$, the PH \eqref{eq:ph} space decomposition is finer
than PAFW: it requires the solution of smaller vector-valued subproblems, in
the stars of edges or faces, instead of in the stars of
vertices~(cf.~\Cref{fig:stars}). For the edge-star problems solved for
$\Hcurl$, in \Cref{fig:spy_feec}(e) we apply a reverse Cuthill--McKee reordering and observe that the matrix
is block-diagonal with $p$ sparse blocks of size $\bigo{p} \times \bigo{p}$. Hence the Cholesky factorization for the PH patch has
optimal storage without requiring the use of incomplete factorizations, assuming fill-in in the entire block.
On the
auxiliary scalar-valued problem posed on the vertex star, we employ the
incomplete Cholesky factorization described above.  For the face-star problems
solved in $\Hdiv$, there is no coupling at all between the face degrees of
freedom in the FDM basis, and ICC(0) offers a direct
solver~(\Cref{fig:spy_feec}(f)). In fact, the SC-PH relaxation for $\Hdiv$ is
equivalent to point-Jacobi applied to the interface Schur complement.

\section{Numerical experiments}
We provide an implementation of the $\P_p$ and $\DP_{p-1}$ elements with the
FDM basis functions on the interval in FIAT \cite{kirby04}.  The extension to
quadrilaterals and hexahedra is achieved by taking tensor-products of the
one-dimensional elements with FInAT \cite{homolya17}.  Code for the
sum-factorized evaluation of the residual is automatically generated by
Firedrake \cite{firedrake, homolya18}, implementing a
Gau\ss--Lobatto quadrature rule with $3(p+1)/2$ points along each direction.
The sparse preconditioner discretizing the auxiliary operator is implemented as
a PETSc~\cite{petsc-user-ref} preconditioner as $\texttt{firedrake.FDMPC}$.

The preconditioners described in \Cref{sec:relaxation} have been presented as
additive.  However, the preconditioners in our experiments are implemented as
hybrid multigrid/Schwarz methods~\cite{lottes05}: they combine patches
additively within each level, and the levels are combined multiplicatively in a
V(1,1)-cycle.  The Cholesky factorization of the patch matrices is computed
using CHOLMOD~\cite{davis08} and the ICC factorization is done with PETSc's own
implementation.
Most of our computations were performed on a single node of the ARCHER2 system,
with two 64-core AMD EPYC 7742 CPUs (2.25 GHz) and 512 GiB of RAM.

Code for all examples has been archived and is available at~\cite{zenodo/pmg-de-rham}.

\subsection{Riesz maps: robust iteration counts and optimal complexity} \label{sec:riesz_maps_experiments}
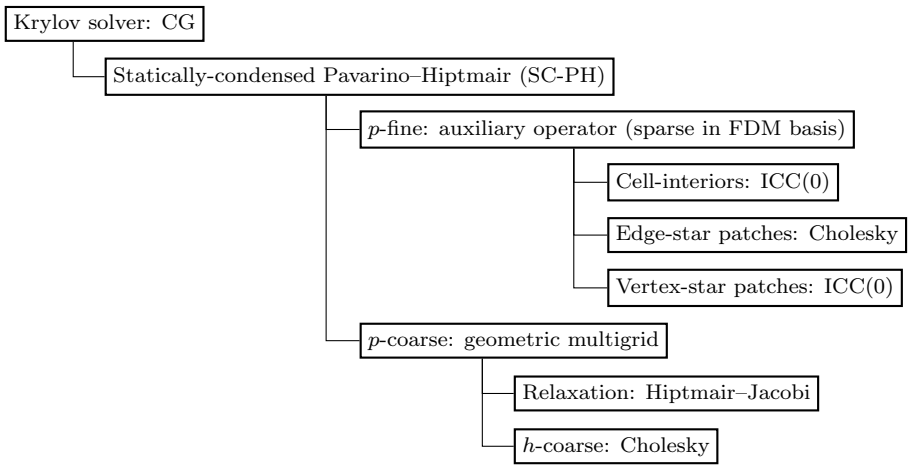
\begin{figure}[tbhp]
\footnotesize
\centering
\begin{tikzpicture}[%
 every node/.style={draw=black, thick, anchor=west},
grow via three points={one child at (-0.0,-0.7) and
	two children at (0.0,-0.7) and (0.0,-1.4)},
edge from parent path={(\tikzparentnode.210) |- (\tikzchildnode.west)}]
\node {Krylov solver: CG}
child {node {Statically-condensed Pavarino--Hiptmair (SC-PH)}
   child {node {$p$-fine: auxiliary operator (sparse in FDM basis)}
      child {node {Cell-interiors: ICC(0)}}
      child {node {Edge-star patches: Cholesky}}
      child {node {Vertex-star patches: ICC(0)}}
   }
   child[missing]{}
   child[missing]{}
   child[missing]{}
   child {node {$p$-coarse: geometric multigrid}
      child {node {Relaxation: Hiptmair--Jacobi}}
      child {node {$h$-coarse: Cholesky}}
   }
};
\end{tikzpicture}
\caption{Solver diagram for the $\Hcurl$ Riesz map
using static condensation and incomplete Cholesky factorization on vertex patches (SC-PH/FDM/ICC).}
\label{fig:hcurl_solver_diagram}
\end{figure}

We first present numerical evidence demonstrating that our preconditioner for the weighted Riesz maps
\eqref{eq:hgrad}--\eqref{eq:hdiv} yields CG iteration counts that are robust to
mesh size $h$, polynomial degree $p$, and the coefficients $\alpha$ and
$\beta$.  For concreteness, in \Cref{fig:hcurl_solver_diagram} we present a
solver diagram for the $H(\curl)$ problem \eqref{eq:hcurl} with the SC-PH
relaxation and geometric multigrid with Hiptmair--Jacobi
relaxation~\cite{hiptmair98} on the $p$-coarse space; the solver diagrams for
the other Riesz maps and space decompositions are analogous.

\begin{figure}[tbhp]
\centering
\hfill
\subfloat[Cartesian]{
   \includegraphics[height=0.3\textwidth,valign=c]{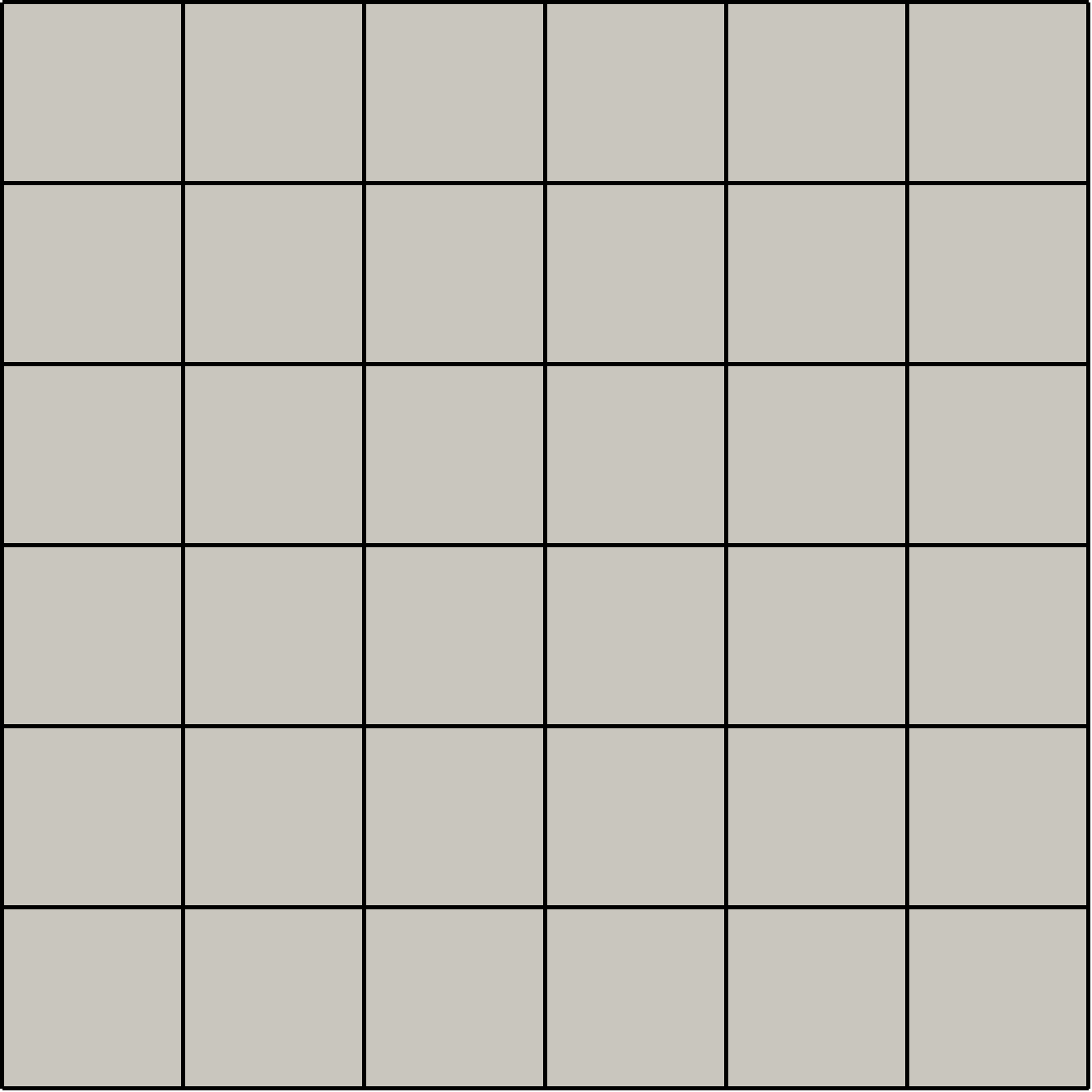}}
\hfill
\subfloat[unstructured]{
   \includegraphics[height=0.3\textwidth,valign=c]{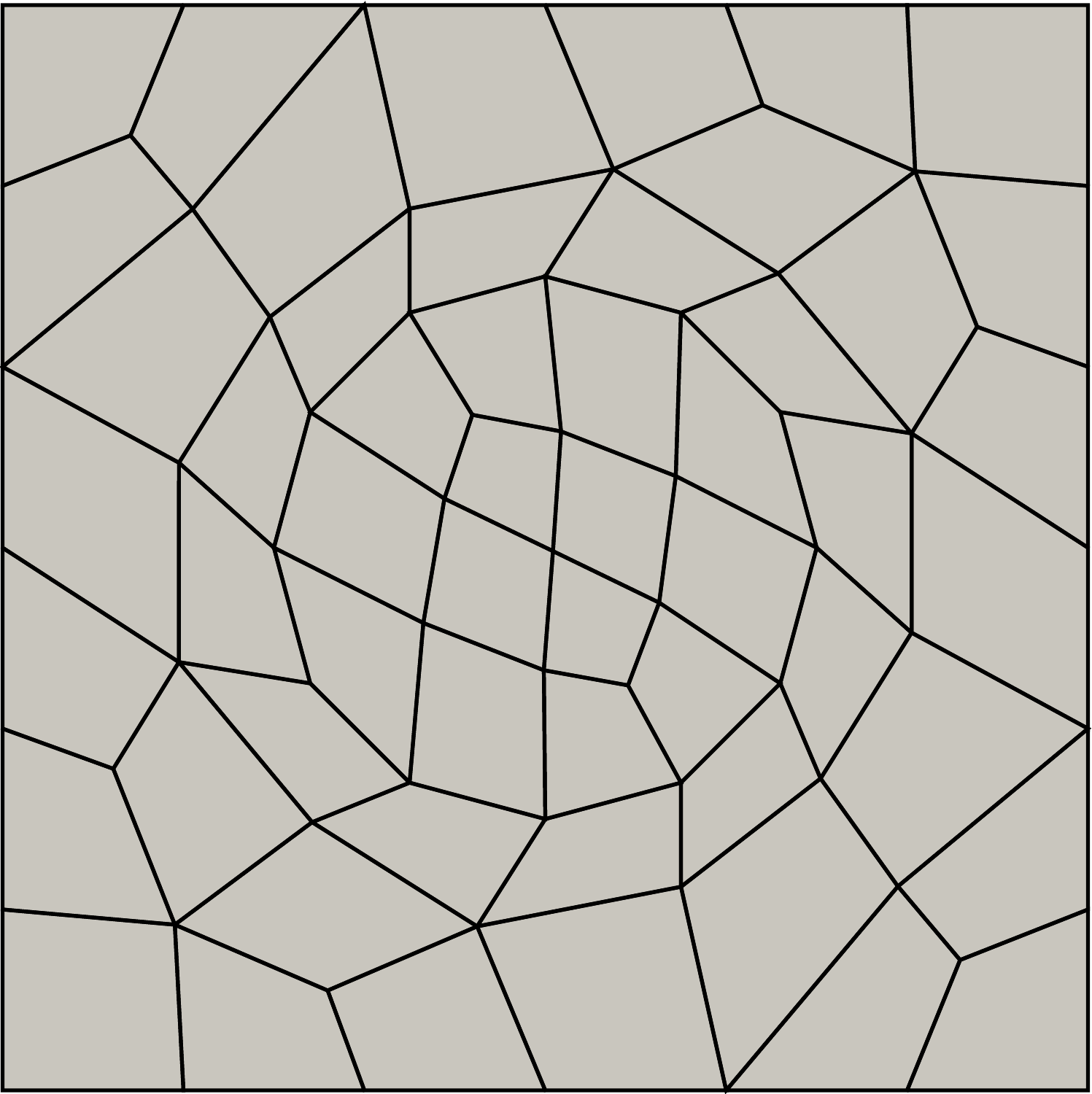}}
\hfill
\caption{The meshes employed in \Cref{sec:riesz_maps_experiments} are the extrusions with six cells in the vertical
of the two-dimensional meshes shown here.}
\label{fig:mesh_hierarchies}
\end{figure}

We consider two mesh hierarchies of $\Omega = (0, 1)^3$, one structured and one
unstructured. 
The coarse meshes are the extrusions with six cells in the
vertical of the two-dimensional meshes depicted in \Cref{fig:mesh_hierarchies}.
We then uniformly refine these $l \ge 0$ times.  Each run is terminated when
the (natural) $P^{-1}$-norm of the residual is reduced by a factor of $10^8$
starting from a zero initial guess. Each problem has homogeneous Dirichlet
boundary conditions on $\Gamma_D = \partial\Omega$ and a randomized right-hand
side that is prescribed independently of $\alpha, \beta$, through its Riesz
representative
\begin{equation} \label{eq:random_rhs}
F(v) = (v, w_{h,p}^k)_\Omega + (\d^k v, \d^k w_{h,p}^k)_\Omega, 
\quad \mbox{for~} w_{h,p}^k \in V^k_{h,p}.
\end{equation}

We first consider the robustness of our solvers with respect to mesh refinement
$l$ and polynomial degree $p$, in \Crefrange{tab:hgrad-iter}{tab:hdiv-iter}. In
these experiments we fix $\alpha = 1$ and $\beta = 10^{-8}$.  The mesh with
distorted cells is useful to measure the effect of the diagonal approximation
in the sparse auxiliary operator \eqref{eq:sum-factorization_approx}, since the
Jacobian of a non-affine coordinate mapping gives rise to variable coefficients
in the reference coordinates, recall \eqref{eq:form-reference-values}.  The results show
almost complete $p$- and $h$-robustness in the Cartesian case, and very slow
growth of iteration counts in the unstructured case.

We next consider the robustness with respect to $\alpha$ and $\beta$ for fixed
$p = 7$, $l = 2$, for the $H(\curl)$ and $H(\div)$ problems
\eqref{eq:hcurl}--\eqref{eq:hdiv}, in
\Crefrange{tab:hcurl-sweep}{tab:hdiv-sweep}. Again, iteration counts do not
vary substantially as the coefficients are varied, indicating the desired
parameter-robustness. As expected, we observe that the iteration counts only
depend on the ratio $\alpha/\beta$, which shows that our solver is invariant to
a rescaling of the system matrix.

\begin{figure}[tbhp]
\centering
 \includegraphics[height=0.3\textwidth,valign=c]{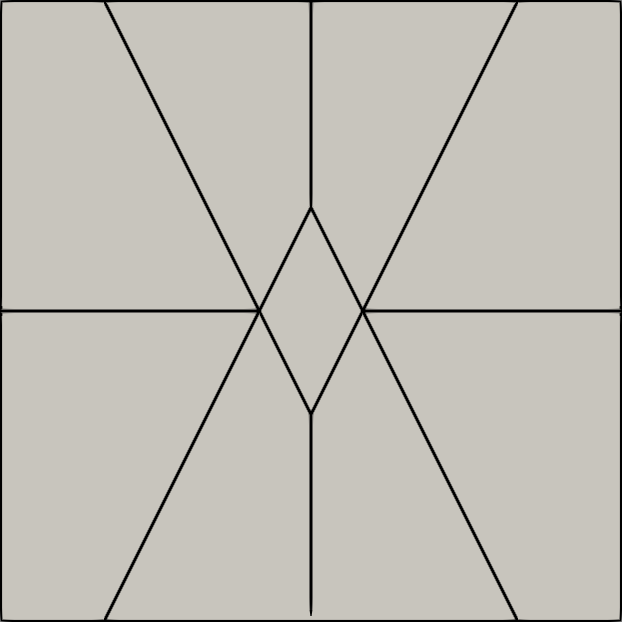}
\caption{
The mesh employed for the complexity plots in 
\Cref{fig:complexities} 
is the extrusion with three cells in the vertical
of the two-dimensional meshes shown here.
}
\label{fig:small_unstructured}
\end{figure}

We next record the flop counts, peak memory usage, and matrix nonzeros for the
PH and SC-PH solvers (with either Cholesky or ICC for the $\Hgrad$ vertex-star
patches) while varying $p$ with $\alpha=\beta=1$ with the mesh shown in
\Cref{fig:small_unstructured}.  This smaller mesh was used so that we could
solve the problem at higher $p$ using the GLL basis. The results were reported
above in \Cref{sec:introduction}, \Cref{fig:complexities}.

\begin{table}[tbhp]
{
\footnotesize
\caption{
CG iteration counts for the $\Hgrad$ Riesz map solved with
the $\{$chol/SC-chol; ICC/SC-ICC$\}$ preconditioners. Empty table cells indicate an unsuccessful termination due to memory limitations.
}
\label{tab:hgrad-iter}
\begin{center}
\rowcolors{4}{}{gray!25}
\setlength{\tabcolsep}{4pt}
\begin{tabular}{r|*{3}{c@{\hskip 4pt}c}|*{3}{c@{\hskip 4pt}c}}
\toprule
& \multicolumn{6}{c|}{Cartesian    } 
& \multicolumn{6}{c}{Unstructured  }
\\
$p\setminus l$ 
& \multicolumn{2}{c}{0} & \multicolumn{2}{c}{1} & \multicolumn{2}{c|}{2} 
& \multicolumn{2}{c}{0} & \multicolumn{2}{c}{1} & \multicolumn{2}{c}{2} 
\\
\midrule
1 & 1/1 & 1/1 & 8/8 & 8/8 & 8/8 & 8/8 & 1/1 & 1/1 & 12/12 & 12/12 & 14/14 & 14/14\\
3 & 14/11 & 14/11 & 15/11 & 15/11 & 15/11 & 15/11 & 18/18 & 18/18 & 18/19 & 18/19 & 18/21 & 18/21\\
7 & 12/9 & 13/10 & 12/9 & 13/10 & 12/9 & 13/10 & 19/21 & 20/21 & 19/23 & 20/23 & 19/24 & 20/24\\
11 & 11/9 & 12/10 & 11/9 & 12/10 & --/9 & --/10 & 20/23 & 22/24 & 19/24 & 22/24&  --/-- & --/25\\
15 & 10/8 & 13/11 & 10/8 & 13/11 &&& 21/24 & 24/25 & 20/25 & 23/25&&\\
\bottomrule
\end{tabular}
\end{center}
}
\end{table}

\begin{table}[tbhp]
{
\footnotesize
\caption{
CG iteration counts for the $\Hcurl$ Riesz map solved with
the $\{$PAFW/SC-PAFW; PH/SC-PH$\}$ preconditioners.
}
\label{tab:hcurl-iter}
\begin{center}
\rowcolors{4}{}{gray!25}
\setlength{\tabcolsep}{4pt}
\begin{tabular}{r|*{3}{c@{\hskip 4pt}c}|*{3}{c@{\hskip 4pt}c}}
\toprule
& \multicolumn{6}{c|}{Cartesian    } 
& \multicolumn{6}{c}{Unstructured  }
\\
$p\setminus l$ 
& \multicolumn{2}{c}{0} & \multicolumn{2}{c}{1} & \multicolumn{2}{c|}{2} 
& \multicolumn{2}{c}{0} & \multicolumn{2}{c}{1} & \multicolumn{2}{c}{2} 
\\
\midrule
1 & 2/2 & 2/2 & 13/13 & 11/11 & 14/14 & 12/12 & 2/2 & 2/2 & 18/18 & 16/16 & 20/20 & 19/19\\
3 & 14/11 & 20/14 & 15/12 & 22/16 & 15/12 & 22/16 & 18/20 & 26/23 & 19/22 & 27/24 & 19/23 & 27/25\\
7 & 13/10 & 20/15 & 14/11 & 21/16 & 14/11 & 21/16 & 21/23 & 33/29 & 21/25 & 32/28 & 21/25 & 32/28\\
11 & 13/10 & 20/15 & 14/11 & 21/16 &&& 23/25 & 36/31 & 22/26 & 34/30 &&\\
15 & 13/10 & 20/16 & 14/12 & 21/17 &&& 24/26 & 37/32 & --/-- & --/31 &&\\
\bottomrule
\end{tabular}
\end{center}
}
\end{table}

\begin{table}[tbhp]
{
\footnotesize
\caption{
CG iteration counts for the $\Hdiv$ Riesz map solved with
the $\{$PAFW/SC-PAFW; PH/SC-PH$\}$ preconditioners.
}
\label{tab:hdiv-iter}
\begin{center}
\rowcolors{4}{}{gray!25}
\setlength{\tabcolsep}{4pt}
\begin{tabular}{r|*{3}{c@{\hskip 4pt}c}|*{3}{c@{\hskip 4pt}c}}
\toprule
& \multicolumn{6}{c|}{Cartesian    } 
& \multicolumn{6}{c}{Unstructured  }
\\
$p\setminus l$ 
& \multicolumn{2}{c}{0} & \multicolumn{2}{c}{1} & \multicolumn{2}{c|}{2} 
& \multicolumn{2}{c}{0} & \multicolumn{2}{c}{1} & \multicolumn{2}{c}{2} 
\\
\midrule
1 & 2/2 & 2/2 & 11/11 & 14/14 & 12/12 & 15/15 & 2/2 & 2/2 & 17/17 & 22/22 & 19/19 & 22/22\\
3 & 9/8 & 15/11 & 9/8 & 16/13 & 10/9 & 16/13 & 15/21 & 24/22 & 16/23 & 25/25 & 18/25 & 26/27\\
7 & 9/7 & 17/13 & 9/8 & 17/13 & 9/9 & 17/14 & 18/25 & 31/27 & 18/26 & 30/28 & 19/27 & 30/29\\
11 & 9/7 & 17/13 & 9/8 & 18/14 &&& 19/27 & 33/29 & 19/28 & 31/29&&\\
15 & 9/8 & 17/13 & 9/8 & 18/14 &&& 19/28 & 34/30 & --/28 & --/30&&\\
\bottomrule
\end{tabular}
\end{center}
}
\end{table}

\begin{table}[tbhp]
{
\footnotesize
\caption{
CG iteration counts for the $\Hcurl$ Riesz map discretized with $p=7$, $l=2$
solved using the $\{$PAFW/SC-PAFW; PH/SC-PH$\}$ hybrid preconditioners.
For the table cell marked --, the PH patch matrix is numerically singular.
}
\label{tab:hcurl-sweep}
\begin{center}
\rowcolors{4}{}{gray!25}
\setlength{\tabcolsep}{4pt}
\begin{tabular}{l|*{3}{c@{\hskip 4pt}c}|*{3}{c@{\hskip 4pt}c}}
\toprule
& \multicolumn{6}{c|}{Cartesian    } 
& \multicolumn{6}{c}{Unstructured  }
\\
$\beta\setminus \alpha$ 
& \multicolumn{2}{c}{$10^{-3}$} & \multicolumn{2}{c}{$10^0$} & \multicolumn{2}{c|}{$10^3$} 
& \multicolumn{2}{c}{$10^{-3}$} & \multicolumn{2}{c}{$10^0$} & \multicolumn{2}{c}{$10^3$} 
\\
\midrule
$10^{-6}$ & 12/9 & 20/15 & 14/11 & 22/16 & 16/24 & 22/18 & 19/27 & 31/28 & 21/27 & 32/28 & 20/25 & --/32\\
$10^{-3}$ & 10/8 & 18/14 & 12/9 & 20/15 & 14/11 & 22/16 & 17/27 & 30/28 & 19/27 & 31/28 & 21/27 & 32/28\\
$10^{ 0}$ & 12/10 & 20/15 & 10/8 & 18/14 & 12/9 & 20/15 & 19/27 & 31/28 & 17/27 & 30/28 & 19/27 & 31/28\\
$10^{ 3}$ & 11/8 & 21/13 & 12/10 & 20/15 & 10/8 & 18/14 & 24/28 & 42/32 & 19/27 & 31/28 & 17/27 & 30/28\\
$10^{ 6}$ & 11/8 & 23/14 & 11/8 & 21/13 & 12/9 & 20/15 & 29/30 & 51/38 & 24/28 & 42/32 & 19/27 & 31/28\\
\bottomrule
\end{tabular}
\end{center}
}
\end{table}

\begin{table}[tbhp]
{
\footnotesize
\caption{
CG iteration counts for the $\Hdiv$ Riesz map discretized with $p=7$, $l=2$
solved using the $\{$PAFW/SC-PAFW; PH/SC-PH$\}$ hybrid preconditioners.
}
\label{tab:hdiv-sweep}
\begin{center}
\rowcolors{4}{}{gray!25}
\setlength{\tabcolsep}{4pt}
\begin{tabular}{l|*{3}{c@{\hskip 4pt}c}|*{3}{c@{\hskip 4pt}c}}
\toprule
& \multicolumn{6}{c|}{Cartesian    } 
& \multicolumn{6}{c}{Unstructured  }
\\
$\beta\setminus \alpha$ 
& \multicolumn{2}{c}{$10^{-3}$} & \multicolumn{2}{c}{$10^0$} & \multicolumn{2}{c|}{$10^3$} 
& \multicolumn{2}{c}{$10^{-3}$} & \multicolumn{2}{c}{$10^0$} & \multicolumn{2}{c}{$10^3$} 
\\
\midrule
$10^{-6}$ & 8/8 & 15/11 & 9/9 & 17/14 & 10/11 & 17/15 & 16/23 & 26/25 & 19/27 & 30/29 & 19/27 & 30/29\\
$10^{-3}$ & 6/6 & 12/9 & 8/8 & 15/11 & 9/9 & 17/14 & 12/18 & 20/19 & 16/23 & 26/25 & 19/27 & 30/29\\
$10^{ 0}$ & 8/7 & 15/11 & 6/6 & 12/9 & 8/8 & 15/11 & 16/22 & 26/23 & 12/18 & 20/19 & 16/23 & 26/25\\
$10^{ 3}$ & 7/5 & 21/11 & 8/7 & 15/11 & 6/6 & 12/9 & 18/28 & 38/28 & 16/22 & 26/23 & 12/18 & 20/19\\
$10^{ 6}$ & 7/4 & 23/10 & 7/5 & 21/11 & 8/7 & 15/11 & 22/31 & 48/33 & 18/28 & 38/28 & 16/22 & 26/23\\
\bottomrule
\end{tabular}
\end{center}
}
\end{table}

\subsection{Piecewise-constant coefficients}
As a test case for our $\Hcurl$ multigrid solver, we consider a definite
Maxwell problem proposed by Kolev \& Vassilevski~\cite[\S 6.2]{kolev09} and
adapted by Pazner et al.~\cite[\S 6.4]{pazner22}.  The problem models
electromagnetic diffusion in an annular copper wire in air, with a
piecewise-constant diffusion coefficient $\beta$ in \eqref{eq:hcurl}.  As in
Pazner et al., we employ a $\Q_3$ coordinate field.  We set $\Gamma_D =
\emptyset$, $\alpha=1$, $\beta_{\text{copper}} = 1$, and vary the diffusion
constant of air $\beta_{\text{air}}$. Since the SC-PH relaxation
exhibits the best performance in the experiments of
\Cref{sec:riesz_maps_experiments}, we only consider this solver here, and as
the $p$-coarse solver we apply a single Hypre AMS~\cite{kolev09} algebraic
multigrid cycle. For the Krylov solver, we set a relative tolerance of
$10^{-8}$. 

We first consider robustness of CG iteration counts to the magnitude of the
jump in the coefficients, in \Cref{tab:maxwell-iter}. The results exhibit
almost perfect robustness across twelve orders of magnitude for
$\beta_{\text{air}}$ and across polynomial degrees between 2 and 14. This
contrasts with \cite[Table 2]{pazner22}, where the iteration counts for the
LOR-AMS solver roughly double from $p=2$ to $p=6$.  We also tabulate the memory
and solve times required as a function of $p$ in \Cref{tab:maxwell-cost} for
fixed $\beta_{\text{air}} = 10^{-6}$.

\begin{table}[tbhp]
{
\footnotesize
\caption{
CG iteration counts for the $\Hcurl$ definite Maxwell problem  
with piecewise-constant coefficients, solved with the SC-PH preconditioner.
}
\label{tab:maxwell-iter}
\begin{center}
\rowcolors{4}{}{gray!25}
\setlength{\tabcolsep}{4pt}
\begin{tabular}{rr|*{5}{c}}
\toprule
& & \multicolumn{5}{c}{$\beta_\text{air}$}\\
$p$ & \#DOFs & $10^{-6}$ & $10^{-3}$ & $10^{0}$ & $10^{3}$ & $10^{6}$ 
\\
\midrule
2  &     516,820 & 25 & 23 & 21 & 29 & 33\\
3  &   1,731,408 & 24 & 22 & 21 & 26 & 31\\
4  &   4,088,888 & 25 & 23 & 21 & 25 & 31\\
5  &   7,968,340 & 26 & 24 & 21 & 28 & 30\\
6  &  13,748,844 & 27 & 24 & 21 & 26 & 29\\
10 &  63,462,980 & 28 & 25 & 21 & 27 & 29\\
14 & 173,920,348 & 28 & 25 & 22 & 28 & 25\\
\bottomrule
\end{tabular}
\end{center}
}
\end{table}

\begin{table}[tbhp]
{
\footnotesize
\caption{
Memory usage and runtimes for the definite Maxwell problem in $\Hcurl$ with
piecewise-constant coefficients ($\beta_\text{air} = 10^{-6}$) solved with
the SC-PH preconditioner.  NNZ Mat includes the number of nonzeros of the
symmetric part of $\Hcurl$ and $\Hgrad$ Schur complements, the ideal
restriction matrix onto the interface of $\Hcurl$, and the tabulation of the
gradient of the interface basis functions.  Memory records the storage of
all matrices across the solver.  The runtime has been broken down into
assembly of the sparse Schur complements, setup of the subspace problems,
and solve times.
}
\label{tab:maxwell-cost}
\begin{center}
\rowcolors{3}{gray!25}{}
\setlength{\tabcolsep}{4pt}
\sisetup{table-auto-round}
\begin{tabular}{r|*{2}{S[table-format=1.2e1]}*{4}{S[table-format=3.2]}}
\toprule
$p$ & 
\multicolumn{1}{c}{NNZ Mat} & 
\multicolumn{1}{c}{NNZ Fact} & 
\multicolumn{1}{c}{Memory (GB)} & 
\multicolumn{1}{c}{Assembly (s)} & 
\multicolumn{1}{c}{Setup (s)} & 
\multicolumn{1}{c}{Solve (s)}  
\\
\midrule
2 & 1.7318e+07 & 1.4738e+07 & 1.03 & 0.42 & 1.31 & 1.11 \\
3 & 5.3952e+07 & 6.5279e+07 & 2.53 & 0.59 & 1.58 & 1.46 \\
4 & 1.2503e+08 & 1.7344e+08 & 5.45 & 0.90 & 3.47 & 2.73 \\
5 & 2.4192e+08 & 3.6076e+08 & 10.22 & 1.42 & 6.46 & 4.48 \\
6 & 4.16e+08 & 6.4875e+08 & 17.28 & 2.22 & 8.49 & 7.64 \\
10 & 1.9117e+09 & 3.2381e+09 & 76.04 & 9.60 & 31.67 & 34.15 \\
14 & 5.2321e+09 & 9.1607e+09 & 209.86 & 27.08 & 75.10 & 94.90 \\
\bottomrule
\end{tabular}
\end{center}
}
\end{table}

\subsection{Mixed formulation of Hodge Laplacians}

The Riesz maps provide building blocks for developing preconditioners for more
complex systems~\cite{hiptmair06,kirby10,mardal11,malek14}. In this final
example, we demonstrate this by constructing preconditioners for the mixed
formulation of the Hodge Laplacians associated with the $\Ltwo$ de Rham
complex. For $\Gamma_D = \emptyset$, the problem is to find $(\sigma, u) \in V^{k-1}\times V^k$ such
that
\begin{alignat}{2}
-(\tau, \sigma)_\Omega + (\d^{k-1}\tau, u)_\Omega &= 0 &&\quad \forall \,\tau \in V^{k-1},\\
(v, \d^{k-1}\sigma)_\Omega + (\d^k v, \d^k u)_\Omega &= F(v) &&\quad \forall \,v \in V^k,
\end{alignat}
where $F(v)$ is a random right-hand side given by~\eqref{eq:random_rhs}.
For $k=3$, the Hodge Laplacian problem corresponds to a mixed formulation of the Poisson equation.

To solve this saddle point problem, we follow the operator preconditioning
framework of Hiptmair and Mardal \& Winther~\cite{hiptmair06,mardal11},
employing a block-diagonal preconditioner with the Riesz maps
($\alpha=\beta=1$) for $V^{k-1}$ and $V^k$.  We use 4 Chebyshev iterations
preconditioned by the SC-PH solver for each block. For $\Ltwo$, we use 4
Chebyshev iterations preconditioned by point-Jacobi in the FDM basis for
$\DQ_{p-1}$. As outer Krylov solver we employ the minimum residual method
(MINRES)~\cite{paige75}, as this allows for the solution of indefinite problems
using a symmetric coercive preconditioner.  The convergence criterion for the
iteration is a relative reduction of the (natural) $P^{-1}$-norm of the
residual by a factor of $10^8$, starting from a zero initial guess. The
iteration counts for the cases $k=1,2,3$ are reported in \Cref{tab:hodge-iter}.
As for the Riesz maps, we observe robustness with respect to both $h$ and $p$.
The 4 inner Chebyshev iterations ensure that the preconditioner for each block
is properly scaled, and for these problems, the MINRES iteration counts are
typically reduced by a factor greater than 4 when compared to those obtained
with a single V-cycle on each block.

\begin{table}[tbhp]
{
\footnotesize
\caption{
MINRES iteration counts for the mixed formulation of 
the Hodge Laplacians $(k=1,2,3)$ preconditioned with the Riesz maps 
using 4 Chebyshev iterations of the SC-PH solvers on each block.
}
\label{tab:hodge-iter}
\begin{center}
\rowcolors{4}{}{gray!25}
\setlength{\tabcolsep}{4pt}
\begin{tabular}{r|*{3}{c}|*{3}{c}||*{3}{c}|*{3}{c}||*{3}{c}|*{3}{c}}
\toprule
& \multicolumn{6}{c||}{$\Hgrad \times \Hcurl$} 
& \multicolumn{6}{c||}{$\Hcurl \times \Hdiv$} 
& \multicolumn{6}{c}{$\Hdiv \times \Ltwo$} 
\\
& \multicolumn{3}{c|}{Cartesian} & \multicolumn{3}{c||}{Unstructured}
& \multicolumn{3}{c|}{Cartesian} & \multicolumn{3}{c||}{Unstructured}
& \multicolumn{3}{c|}{Cartesian} & \multicolumn{3}{c}{Unstructured}
\\
$p\setminus l$ 
& \multicolumn{1}{c}{0} & \multicolumn{1}{c}{1} & \multicolumn{1}{c|}{2} 
& \multicolumn{1}{c}{0} & \multicolumn{1}{c}{1} & \multicolumn{1}{c||}{2} 
& \multicolumn{1}{c}{0} & \multicolumn{1}{c}{1} & \multicolumn{1}{c|}{2} 
& \multicolumn{1}{c}{0} & \multicolumn{1}{c}{1} & \multicolumn{1}{c||}{2} 
& \multicolumn{1}{c}{0} & \multicolumn{1}{c}{1} & \multicolumn{1}{c|}{2} 
& \multicolumn{1}{c}{0} & \multicolumn{1}{c}{1} & \multicolumn{1}{c}{2} 
\\
\midrule
3  & 9 & 9 & 9 & 16 & 16 & 17 & 8 & 7 & 7 & 15 & 16 & 18 & 8 & 8 & 8 & 12 & 14 & 14\\
7  & 9 & 9 & 8 & 20 & 19 & 19 & 8 & 8 & 6 & 19 & 18 & 18 & 8 & 8 & 8 & 14 & 14 & 13\\
11 & 9 & 9 &   & 23 & 21 &    & 8 & 8 &   & 20 & 18 &    & 8 & 8 &   & 15 & 14 &   \\  
15 & 9 & 9 &   & 23 & 21 &    & 8 & 6 &   & 20 & 18 &    & 8 & 8 &   & 15 & 14 &   \\ 
\bottomrule
\end{tabular}
\end{center}
}
\end{table}

\section{Conclusion} \label{sec:conclusion}

We have developed multigrid solvers for the Riesz maps associated with the
$L^2$ de Rham complex whose space and time complexities in polynomial degree
are the same as that required for operator application. Numerical experiments
demonstrate that the solvers are robust to mesh refinement, polynomial degree,
and problem coefficients, and that they remain effective on unstructured grids.
However, the solvers are not robust with respect to anisotropy, in common
with other methods \cite{ceed21}.
The approach relies on developing new finite elements with desirable
interior-orthogonality properties, auxiliary operators that are sparse by
construction, the careful use of incomplete factorizations, and the choice of
space decomposition. The resulting solvers can be employed in the operator
preconditioning framework to develop preconditioners for more complex problems
with solution variables in $\Hgrad$, $\Hcurl$, and $\Hdiv$.

\bibliographystyle{siamplain}
\bibliography{references}

\begin{thebibliography}{10}

\bibitem{anderson99}
{\sc E.~Anderson, Z.~Bai, C.~Bischof, S.~Blackford, J.~Demmel, J.~Dongarra,
  J.~Du~Croz, A.~Greenbaum, S.~Hammarling, A.~McKenney, and D.~Sorensen}, {\em
  {LAPACK} Users' Guide}, SIAM, Philadelphia, PA, third~ed., 1999.

\bibitem{arnold18}
{\sc D.~N. Arnold}, {\em {Finite Element Exterior Calculus}}, SIAM, 2018.

\bibitem{arnold15}
{\sc D.~N. Arnold, D.~Boffi, and F.~Bonizzoni}, {\em Finite element
  differential forms on curvilinear cubic meshes and their approximation
  properties}, Numer. Math., 129 (2015), pp.~1--20.

\bibitem{arnold97}
{\sc D.~N. Arnold, R.~Falk, and R.~Winther}, {\em Preconditioning in
  {$H(\text{div})$} and applications}, Math. Comput., 66 (1997), pp.~957--984.

\bibitem{arnold00}
{\sc D.~N. Arnold, R.~S. Falk, and R.~Winther}, {\em {Multigrid in
  $H(\mathrm{div})$ and $H(\mathrm{curl})$}}, Numer. Math., 85 (2000),
  pp.~197--217.

\bibitem{petsc-user-ref}
{\sc S.~Balay, S.~Abhyankar, M.~F. Adams, S.~Benson, J.~Brown, P.~Brune,
  K.~Buschelman, E.~Constantinescu, L.~Dalcin, A.~Dener, V.~Eijkhout,
  J.~Faibussowitsch, W.~D. Gropp, V.~Hapla, T.~Isaac, P.~Jolivet, D.~Karpeev,
  D.~Kaushik, M.~G. Knepley, F.~Kong, S.~Kruger, D.~A. May, L.~C. McInnes,
  R.~T. Mills, L.~Mitchell, T.~Munson, J.~E. Roman, K.~Rupp, P.~Sanan,
  J.~Sarich, B.~F. Smith, S.~Zampini, H.~Zhang, H.~Zhang, and J.~Zhang}, {\em
  {PETSc/TAO} users manual}, Tech. Report ANL-21/39 - Revision 3.19, Argonne
  National Laboratory, 2023.

\bibitem{brenner08}
{\sc S.~C. Brenner and L.~R. Scott}, {\em The Mathematical Theory of Finite
  Element Methods}, vol.~15 of Texts in Applied Mathematics, Springer-Verlag
  New York, third~ed., 2008.

\bibitem{brubeck22}
{\sc P.~D. Brubeck and P.~E. Farrell}, {\em A scalable and robust vertex-star
  relaxation for high-order {FEM}}, SIAM J. Sci. Comput., 44 (2022),
  pp.~A2991--A3017.

\bibitem{canuto10}
{\sc C.~Canuto, P.~Gervasio, and A.~Quarteroni}, {\em Finite-element
  preconditioning of {G-NI} spectral methods}, SIAM J. Sci. Comput., 31 (2010),
  pp.~4422--4451.

\bibitem{davis08}
{\sc Y.~Chen, T.~A. Davis, W.~W. Hager, and S.~Rajamanickam}, {\em {Algorithm
  887: CHOLMOD, Supernodal Sparse Cholesky Factorization and Update/Downdate}},
  ACM Trans. Math. Softw., 35 (2008).

\bibitem{ciarlet02}
{\sc P.~G. Ciarlet}, {\em {The finite element method for elliptic problems}},
  SIAM, 2002.

\bibitem{deville90}
{\sc M.~O. Deville and E.~H. Mund}, {\em Finite-element preconditioning for
  pseudospectral solutions of elliptic problems}, SIAM J. Sci. Stat. Comput.,
  11 (1990), pp.~311--342.

\bibitem{dohrmann21}
{\sc C.~R. Dohrmann}, {\em {Spectral Equivalence of Low-Order Discretizations
  for High-Order $H(\mathrm{curl})$ and $H(\mathrm{div})$ Spaces}}, SIAM J.
  Sci. Comput., 43 (2021), pp.~A3992--A4014.

\bibitem{egorova09}
{\sc O.~Egorova, M.~Savchenko, V.~Savchenko, and I.~Hagiwara}, {\em Topology
  and geometry of hexahedral complex: combined approach for hexahedral
  meshing}, J. Comput. Sci. Technol., 3 (2009), pp.~171--182.

\bibitem{falgout02}
{\sc R.~D. Falgout and U.~M. Yang}, {\em {Hypre: A Library of High Performance
  Preconditioners}}, in Computational Science -- ICCS 2002, P.~M.~A. Sloot,
  A.~G. Hoekstra, C.~J.~K. Tan, and J.~J. Dongarra, eds., vol.~2331 of Lecture
  Notes in Computer Science, Springer Berlin Heidelberg, 2002, pp.~632--641.

\bibitem{pcpatch}
{\sc P.~E. Farrell, M.~G. Knepley, L.~Mitchell, and F.~Wechsung}, {\em
  {PCPATCH}: software for the topological construction of multigrid relaxation
  methods}, ACM Trans. Math. Softw., 47 (2021).

\bibitem{firedrake}
{\sc D.~A. Ham, P.~H.~J. Kelly, L.~Mitchell, C.~J. Cotter, R.~C. Kirby,
  K.~Sagiyama, N.~Bouziani, S.~Vorderwuelbecke, T.~J. Gregory, J.~Betteridge,
  D.~R. Shapero, R.~W. Nixon-Hill, C.~J. Ward, P.~E. Farrell, P.~D. Brubeck,
  I.~Marsden, T.~H. Gibson, M.~Homolya, T.~Sun, A.~T.~T. McRae, F.~Luporini,
  A.~Gregory, M.~Lange, S.~W. Funke, F.~Rathgeber, G.-T. Bercea, and G.~R.
  Markall}, {\em Firedrake user manual},  (2023).

\bibitem{hientzsch01}
{\sc B.~Hientzsch}, {\em Fast solvers and domain decomposition preconditioners
  for spectral element discretizations of problems in $H(\mathrm{curl})$}, PhD
  thesis, New York University, 2001.

\bibitem{hientzsch05}
{\sc B.~Hientzsch}, {\em Domain decomposition preconditioners for spectral
  {N}{\'e}d{\'e}lec elements in two and three dimensions}, in Domain
  Decomposition Methods in Science and Engineering, Springer, 2005,
  pp.~597--604.

\bibitem{hiptmair98}
{\sc R.~Hiptmair}, {\em {Multigrid Method for Maxwell's Equations}}, SIAM J.
  Numer. Anal., 36 (1998), pp.~204--225.

\bibitem{hiptmair06}
{\sc R.~Hiptmair}, {\em Operator preconditioning}, Comput. Math. Appl., 52
  (2006), pp.~699--706.

\bibitem{hiptmair07}
{\sc R.~Hiptmair and J.~Xu}, {\em {Nodal auxiliary space preconditioning in
  $H(\mathrm{curl})$ and $H(\mathrm{div})$ spaces}}, SIAM J. Numer. Anal., 45
  (2007), pp.~2483--2509.

\bibitem{homolya17}
{\sc M.~Homolya, R.~C. Kirby, and D.~A. Ham}, {\em Exposing and exploiting
  structure: optimal code generation for high-order finite element methods},
  2017, \url{https://arxiv.org/abs/1711.02473}.

\bibitem{homolya18}
{\sc M.~Homolya, L.~Mitchell, F.~Luporini, and D.~A. Ham}, {\em {TSFC}: a
  structure-preserving form compiler}, SIAM J. Sci. Comput., 40 (2018),
  pp.~C401--C428.

\bibitem{kirby04}
{\sc R.~C. Kirby}, {\em {Algorithm 839: FIAT, a new paradigm for computing
  finite element basis functions}}, ACM Trans. Math. Soft., 30 (2004),
  pp.~502--516.

\bibitem{kirby10}
{\sc R.~C. Kirby}, {\em From functional analysis to iterative methods}, SIAM
  Rev., 52 (2010), pp.~269--293.

\bibitem{knepley09}
{\sc M.~G. Knepley and D.~A. Karpeev}, {\em {Mesh algorithms for PDE with Sieve
  I: Mesh distribution}}, Sci. Program., 17 (2009), pp.~215--230.

\bibitem{ceed21}
{\sc T.~Kolev et~al.}, {\em {High-order algorithmic developments and
  optimizations for large-scale GPU-accelerated simulations}}, Tech. Report
  {CEED-MS36}, Lawrence Livermore National Laboratory, Livermore, CA, 2021.

\bibitem{kolev21}
{\sc T.~Kolev, P.~Fischer, M.~Min, J.~Dongarra, J.~Brown, V.~Dobrev,
  T.~Warburton, S.~Tomov, M.~S. Shephard, A.~Abdelfattah, et~al.}, {\em
  Efficient exascale discretizations: High-order finite element methods}, Int.
  J. High Perform. Comput. Appl., 35 (2021), pp.~527--552.

\bibitem{kolev09}
{\sc T.~V. Kolev and P.~S. Vassilevski}, {\em Parallel auxiliary space {AMG}
  for {$H(\mathrm{curl})$} problems}, J. Comput. Math., 27 (2009),
  pp.~604--623.

\bibitem{kolev12}
{\sc T.~V. Kolev and P.~S. Vassilevski}, {\em Parallel auxiliary space {AMG}
  solver for {$H(\mathrm{div})$} problems}, SIAM J. Sci. Comput., 34 (2012),
  pp.~A3079--A3098.

\bibitem{lee07}
{\sc Y.-J. Lee, J.~Wu, J.~Xu, and L.~Zikatanov}, {\em {Robust} {subspace}
  {correction} {methods} {for} {nearly} {singular} {systems}}, Math. Models
  Methods Appl. Sci., 17 (2007), pp.~1937--1963.

\bibitem{logg09}
{\sc A.~Logg}, {\em Efficient representation of computational meshes}, Int. J.
  Comput. Sci. Eng., 4 (2009), pp.~283--295.

\bibitem{lottes05}
{\sc J.~W. Lottes and P.~F. Fischer}, {\em {Hybrid multigrid/Schwarz algorithms
  for the spectral element method}}, J. Sci. Comput., 24 (2005), pp.~45--78.

\bibitem{lynch64}
{\sc R.~E. Lynch, J.~R. Rice, and D.~H. Thomas}, {\em Direct solution of
  partial difference equations by tensor product methods}, Numer. Math., 6
  (1964), pp.~185--199.

\bibitem{malek14}
{\sc J.~M\'alek and Z.~Strako\v{s}}, {\em Preconditioning and the Conjugate
  Gradient Method in the Context of Solving PDEs}, vol.~1 of SIAM Spotlights,
  SIAM, 2014.

\bibitem{mardal11}
{\sc K.-A. Mardal and R.~Winther}, {\em Preconditioning discretizations of
  systems of partial differential equations}, Numer. Linear Algebra Appl., 18
  (2011), pp.~1--40.

\bibitem{munkres84}
{\sc J.~R. Munkres}, {\em Elements of Algebraic Topology}, {CRC} Press, 1984.

\bibitem{nedelec80}
{\sc J.-C. N{\'e}d{\'e}lec}, {\em Mixed finite elements in $\mathbb{R}^3$},
  Numer. Math., 35 (1980), pp.~315--341.

\bibitem{orszag80}
{\sc S.~A. Orszag}, {\em Spectral methods for problems in complex geometries},
  J. Comput. Phys., 37 (1980), pp.~70--92.

\bibitem{paige75}
{\sc C.~C. Paige and M.~A. Saunders}, {\em Solution of sparse indefinite
  systems of linear equations}, {SIAM} J. Numer. Anal., 12 (1975),
  pp.~617--629.

\bibitem{pavarino93}
{\sc L.~F. Pavarino}, {\em Additive {S}chwarz methods for the $p$-version
  finite element method}, Numer. Math., 66 (1993), pp.~493--515.

\bibitem{pazner20}
{\sc W.~Pazner}, {\em Efficient low-order refined preconditioners for
  high-order matrix-free continuous and discontinuous {G}alerkin methods}, SIAM
  J. Sci. Comput., 42 (2020), pp.~A3055--A3083.

\bibitem{pazner22}
{\sc W.~Pazner, T.~Kolev, and C.~R. Dohrmann}, {\em {Low-Order Preconditioning
  for the High-Order Finite Element de Rham Complex}}, 45 (2023),
  pp.~A675--A702.

\bibitem{pellikka13}
{\sc M.~Pellikka, S.~Suuriniemi, L.~Kettunen, and C.~Geuzaine}, {\em Homology
  and cohomology computation in finite element modeling}, SIAM J. Sci. Comput.,
  35 (2013), pp.~B1195--B1214.

\bibitem{phillips23}
{\sc M.~Phillips and P.~Fischer}, {\em {Optimal Chebyshev Smoothers and
  One-sided V-cycles}}, 2023, \url{https://arxiv.org/abs/2210.03179}.

\bibitem{schoberl99}
{\sc J.~Sch{\"o}berl}, {\em Robust multigrid methods for parameter dependent
  problems}, PhD thesis, Johannes Kepler Universität Linz, Linz, Austria,
  1999.

\bibitem{schoberl08}
{\sc J.~Sch{\"o}berl, J.~M. Melenk, C.~Pechstein, and S.~Zaglmayr}, {\em
  Additive {S}chwarz preconditioning for p-version triangular and tetrahedral
  finite elements}, IMA J. Numer. Anal., 28 (2008), pp.~1--24.

\bibitem{schoberl05}
{\sc J.~Sch{\"o}berl and S.~Zaglmayr}, {\em High order {N}{\'e}d{\'e}lec
  elements with local complete sequence properties}, COMPEL - Int. J. Comput.
  Math. Electr. Electron. Eng.,  (2005).

\bibitem{schwedes17}
{\sc T.~Schwedes, D.~A. Ham, S.~W. Funke, and M.~D. Piggott}, {\em Mesh
  Dependence in {PDE}-Constrained Optimisation}, Springer, 2017.

\bibitem{thompson23}
{\sc J.~L. Thompson, J.~Brown, and Y.~He}, {\em {Local Fourier Analysis of
  p-Multigrid for High-Order Finite Element Operators}}, 45 (2023),
  pp.~S351--S370.

\bibitem{xu92}
{\sc J.~Xu}, {\em Iterative methods by space decomposition and subspace
  correction}, SIAM Rev., 34 (1992), pp.~581--613.

\bibitem{zaglmayr06}
{\sc S.~Zaglmayr}, {\em High Order Finite Element Methods for Electromagnetic
  Field Computation.}, PhD thesis, Johannes Kepler University Linz, 2006.

\bibitem{zenodo/pmg-de-rham}
{\em {Software used in `Multigrid solvers for the de Rham complex with optimal
  complexity in polynomial degree'}}.
\newblock \url{https://doi.org/10.5281/zenodo.7358044}, Nov 2022.

\end{thebibliography}
\end{document}